\documentclass[twoside,11pt,reqno]{amsart}

\usepackage{amsmath,amsthm,amssymb,amstext,amsfonts,amscd}
\usepackage{graphicx}
\usepackage{multirow}
\usepackage{url}

\newtheorem{theorem}{Theorem}
\newtheorem{definition}{Definition}
\newtheorem{corollary}{Corollary}
\newtheorem{remark}{Remark}

\numberwithin{equation}{section}

\begin{document}
\title[{\bf On extension of the $\;_{r}R_{s,q}(\alpha,\beta,z)$ function and their $q$-calculus}]
{On extension of the $\;_{r}R_{s,q}(\alpha,\beta,z)$ function and their $q$-calculus}

\author[{\bf A. Shehata and D. Kumar}]{\bf Ayman Shehata and Dinesh Kumar$^{*}$}

\address{A. Shehata: Department of Mathematics, Faculty of Science, Assiut University, Assiut--71516,  EGYPT}
\email{drshehata2006@yahoo.com; aymanshehata@science.aun.edu.eg}

\address{D. Kumar: Department of Applied Sciences, College of Agriculture--Jodhpur,  Agriculture University Jodhpur, Jodhpur--342304, INDIA}
\email{dinesh\_dino03@yahoo.com}

\bigskip

\thanks{$^*$ Corresponding author}



\begin{abstract}
In this article, we investigate and establish some properties including analytic properties, contiguous relations, differential properties, differential operators, an expansion formula, and simple integrals, integral operators, some fractional integral properties, some new integral representations, the Riemann–Liouville (R–L) fractional $q$-derivative and $q$-integral operators of $q$-analogue of various basic $\;_{r}R_{s,q}$ function by using technique of $q$-calculus. Certain interesting consequences of the theorem are also discussed by considering some examples.
\end{abstract}
\maketitle

\section {Introduction and Preliminaries}\label{sec 1}
The theory of basic or $q$-hypergeometric functions and fractional calculus has a wide range of applications in various fields of Mathematical analysis, Engineering, Physical and Sciences, namely-Number theory, Partition theory, Combinatorial analysis, Lie theory, Quantum theory etc. Abdeljawad and Baleanu \cite{ab} examined initial value problems involving the Caputo $q$-fractional derivatives and introduced a $q$-analog of the Mittag-Leffler function. In a separate work, Abdeljawad et al. \cite{abb} explored a generalized form of the $q$-Mittag-Leffler function in the context of $q$-Caputo fractional linear equations. Agarwal \cite{ag} investigated specific fractional $q$-integrals and $q$-derivatives. Ali and Suthar \cite{as} proposed a Riemann–Liouville fractional $q$-calculus operator that incorporates the $q$-Mittag–Leffler function. Kumar et al. \cite{KumarAyantUcarPurohit2022} gave the
certain fractional $q$--integral formulas for the basic $I$--function of two variables.
Sahni et al. \cite{SahniKumarAyantSingh2021} studied a transformation involving basic multivariable $I$--function of Prathima. Furthermore, Kumar et al. \cite{KumarBaleanuAyantSudland2022} studied transformations involving the basic analogue of the Aleph-function of two variables. Ben Said and El-Kamel \cite{be} established the product formula for the generalized $q$-Bessel function. Desai and Shukla \cite{ds} introduced Some results on function $\;_{r}R_{s,q}(a,b,z)$. Heragy et al. \cite{hmo}  evaluated Indefinite $q$-integrals from a method using $q$-Riccati equations. Luca \cite{lu} investigated the positive solutions of a system of fractional $q$-difference equations subject to multi-point boundary conditions. Ganie et al. \cite{gar} established the certain $q$-integrals involving the generalized hypergeometric and basic hypergeometric functions. Swarttouw \cite{sa}  discussed the contiguous function relations for basic hypergeometric series. Shimelis and Suthar \cite{ss1} described On Saigo fractional $q$-Calculus of a general class of $q $-Polynomials. Recently, Kumar et al. \cite{KumarAyantTariboon2020} studied transformations involving the basic analogue of the multivariable $H$-function. Further, Kumar et al. \cite{KumarAyantSudlandChoi2023} studied fractional $q$-integral operators involving the basic analogue of the multivariable Aleph-function.
 Mansour \cite{ma} demonstrated the Linear sequential $q$-difference equations of fractional order. Gonzalez et al. \cite{gms} evaluated Some $q$-generating functions for associated generalized hypergeometric polynomials. Recently, Cai et al. \cite{ckms} studied some properties of the bibasic Humbert hypergeometric functions $E_{1}$ and $E_{2}$ in \cite{ckms}. Early, Shahata discussed the basic Horn hypergeometric functions $H_{1}$, $H_{4}$, $H_{6}$ and $H_{7}$ in \cite{s1, s2}.

In the theory of $q$-series (see \cite{fi,gr,am,ha}), for a real or complex numbers $\alpha$ and $0<|q|<1$, $q\in \mathbb{C}$, the $q$-shifted factorial is defined by
\begin{equation} \label{1.1}
\begin{split}
\left(\alpha;q\right)_{\iota}=\left\{
\begin{array}{ll}
1, & \hbox{$n=0$} \\
\Pi_{i=0}^{\iota-1}(1-aq^{i}), & \hbox{$\iota\in \mathbb{N}$}
\end{array}
\right.
, \end{split}
\end{equation}
where $\mathbb{N}$ is a set of positive integer.
\begin{definition} \label{1.2}
The $q$-analogue of the gamma function is defined as
\begin{equation}
\begin{split}
\Gamma_{q}(n)=\int_{0}^{\frac{1}{1-q}}t^{n-1}E_{q}(-qt)d_{q}t=\int_{0}^{\frac{\infty}{1-q}}t^{n-1}E_{q}(-qt)\mathrm{d}_{q}t\;\; \left(\Re(n)>0\right),
\end{split}
\end{equation}
where
\begin{equation*}
\begin{split}
E_{q}(t)=\sum_{n=0}^{\infty}\frac{q^{\frac{n(n-1)}{2}}}{[n]_{q}!}t^{n}.
\end{split}
\end{equation*}
\end{definition}
\begin{definition}
The $q$-analogue of $q$-beta function is defined by
\begin{equation} \label{1.3}
\begin{split}
B_{q}\left(a,b\right)=&\frac{\Gamma_{q}(a)\Gamma_{q}(b)}{\Gamma_{q}(a+b)}=\int_{0}^{1}t^{a-1}\left(1-qt\right)_{q}^{b-1}\mathrm{d}_{q}t\\
=&\int_{0}^{1}t^{a-1}\frac{\left(qt;q\right)_{\infty}}{\left(q^{b}t;q\right)_{\infty}}\mathrm{d}_{q}t=\int_{0}^{1}t^{a-1}(qt;q)_{b-1}\mathrm{d}_{q}t\;\;\; \left(\Re(a)>0,\Re(b)>0\right).
	\end{split}
\end{equation}
\end{definition}
\begin{definition}
The $q$-derivative of a function $\Omega(z)$ defined as (see, \cite{ag})
\begin{equation} \label{1.4}
\begin{split}
\left(\mathbf{D}_{q}\Omega\right)(z)=\frac{\Omega(z)-\Omega(qz)}{(1-q)z}\;\;\; \left(z\neq 0\right).
\end{split}
\end{equation}
\end{definition}
\begin{definition}
The fractional $q$-integral and $q$-differential operators of order $\alpha$ defined as (see \cite{ag, am, kst})
\begin{equation} \label{1.5}
\begin{split}
\mathbf{I}_{q}^{\alpha}f(z)=\frac{1}{\Gamma_{q}(\alpha)}\int_{0}^{z}\left(z-qt\right)^{\alpha-1}f(t)\mathrm{d}_{q}t\;\;\; \left(\Re(\alpha)>0\right),
\end{split}
\end{equation}
and
\begin{equation}  \label{1.6}
\begin{split}
\mathbf{D}_{q}^{\alpha}f(z)=D_{q}^{\mu}\left[\mathbf{I}_{q}^{\mu-\alpha}f(z)\right]\;\;\; \left(\Re(\alpha)>0\right),
\end{split}
\end{equation}
	where $\mu$ is the smallest integer such that $\mu>\Re(\alpha)$.
\end{definition}
\begin{definition}
Let $0<\alpha\leq 1$. The left- and right-sided Riemann-Liouville fractional $q$-integral operators \cite{am,as} are defined by the following formulas:	
\begin{equation} \label{1.7}
\mathbf{I}_{b^{-},q}^{\alpha}f(z)=\frac{1}{\Gamma_{q}(\alpha)}\int_{qz}^{b}t^{\alpha-1}\left(\frac{zq}{t};q\right)_{\alpha-1}f(t)\mathrm{d}_{q}t,
\end{equation}
and
\begin{equation} \label{1.8}
\mathbf{I}_{a^{+},q}^{\alpha}f(z)=\frac{z^{\alpha-1}}{\Gamma_{q}(\alpha)}\int_{a}^{z}\left(\frac{tq}{z};q\right)_{\alpha-1}f(t)\mathrm{d}_{q}t.
\end{equation}
\end{definition}
\begin{definition}
For $\alpha>0$ and $[\alpha]=n$, the left- and right-sided Riemann-Liouville fractional $q$-derivatives of order $\alpha$ \cite{am, as, gar}, are given by
\begin{equation} \label{1.9}
\mathbf{D}_{b^{-},q}^{\alpha}f(z)=\left(-\frac{1}{q}\right)^{n}D^{n}_{q^{-1}}\,\mathbf{I}_{b^{-},q}^{n-\alpha}f(z),
\end{equation}
and
\begin{equation} \label{1.10}
\mathbf{D}_{a^{+},q}^{\alpha}f(z)=D^{n}_{q}\,\mathbf{I}_{a^{+},q}^{n-\alpha}f(z).
\end{equation}
\end{definition}
\begin{definition}
For finite positive integers $r$ and $s$,the basic hypergeometric function is defined as (see, \cite{ds,gar,sa})
\begin{equation} \label{1.11}
\begin{split}
&\;_{r}\phi_{s}(a_{1},a_{2},\ldots,a_{r};b_{1},b_{2},\ldots,b_{s};z)\\
&	=\sum_{\iota=0}^{\infty}\frac{z^{\iota}}{\left(q;q\right)_{\iota}}\frac{\prod_{i=1}^{r}(a_{i};q)_{\iota}}{\prod_{j=1}^{q}(b_{j};s)_{\iota}}\left[(-1)^{\iota}q^{\binom{\iota}{2}}\right]^{1+s-r},
\end{split}
\end{equation}
where $a_{i}\; (1\leq i\leq r)$ and $b_{j}\; (1\leq j\leq s)$ are complex numbers.
\begin{enumerate}
\item If $0<|q|<1$ and $r\leq s+1$, the series $\;_{r}\phi_{s}$ \eqref{1.11} converges absolutely for all $|z|<1$,
\item If $r\leq s$, and for any value of $z$,
\item This series (\ref{1.11}) converges absolutely when $|q|>1$ and $|z|<\frac{\left|b_{1}b_{2}\ldots b_{s}q\right|}{\left|a_{1}a_{2} \ldots a_{r}\right|}$,
\item The series diverges for $z\neq0$ when $0<|q|<1$ and $r>s+1$ as well as when $\left|q\right|>1$ and $\left|z\right|>\frac{\left|b_{1}b_{2}\ldots b_{s}q\right|}{\left|a_{1}a_{2}\ldots a_{r}\right|}$, unless it terminates.
\end{enumerate}
\end{definition}

Mansour \cite{ma} introduced a novel form of the $q$-analogue of the $q$-Mittag-Leffler function, defined as
\begin{equation} \label{1.12}
E_{\alpha,\beta}\left(z;q\right)=\sum_{i=0}^{\infty}\frac{1}{\Gamma_{q}\left(\alpha i+\beta\right)}z^{i}\;\;\; \left(\left|z\right|<(1-q)^{-\alpha}\right),
\end{equation}
where $\alpha>0$, $\beta \in \mathbb{C}$.
Recently, Sharma and Jain \cite{sj} presented the $q$-analogue of the generalized $q$-Mittag-Leffler function, given by
\begin{equation} \label{1.13}
E_{\alpha,\beta}^{\eta}(z;q)=\sum_{i=0}^{\infty}\frac{\left(q^{\eta};q\right)_{i}}{\left(q;q\right)_{i}\Gamma_{q}\left(\alpha i+\beta\right)}z^{i},
\end{equation}
where $\alpha,\beta,\eta \in \mathbb{C}$, $\Re(\alpha)>0$, $\Re(\beta)>0$ and $\Re(\eta)>0$.\\
Subsequently, Nadeem \cite{Nadeemetal2020} studied the generalized form of the $q$-Mittag-Leffler function, defined as
\begin{equation}
\begin{split}
E_{\alpha,\beta}^{\delta,\eta}(z;q)=\sum_{i=0}^{\infty}\frac{B_{q}(\delta+i,\eta-\delta)}{B_{q}(\delta,\eta-\delta)}\frac{(q^{\eta};q)_{i}}{(q;q)_{i}\Gamma_{q}(\alpha i+\beta)}z^{i},\label{1.14}
\end{split}
\end{equation}
where $\alpha,\beta,\eta \in \mathbb{C}$, $\Re(\alpha)>0$, $\Re(\beta)>0$, $\Re(\eta)>0$ and $\Re(\eta)>\Re(\delta)>0$.
\begin{definition}
The $q$-Laplace transform of an appropriate function is expressed using the following $q$-integral (see, \cite{nh}):
\begin{equation} \label{1.15}
L_{s,q}\left[f(t)\right]=\frac{1}{1-q}\int_{0}^{s^{-1}}E_{q}^{qst}f(t)\,\mathrm{d}_{q}t,
\end{equation}
and
\begin{equation} \label{1.16}
\ell_{s,q}\left[f(t)\right]=\frac{1}{1-q}\int_{0}^{\infty}e_{q}^{-st}f(t)\,\mathrm{d}_{q}t.
\end{equation}
\end{definition}
\begin{definition}
The $q$-extension of the exponential functions are defined by (see \cite{ss})
\begin{equation} \label{1.17}
E_{q}^{z}=\sum_{i=0}^{\infty}\frac{q^{\binom{i}{2}}}{\left(q;q\right)_{i}}z^{i}=\left(-z;q\right)_{\infty}=\;_{0}\Phi_{0}\left(-;-;q,-z\right),
\end{equation}
and
\begin{equation} \label{1.18}	e_{q}^{z}=\sum_{\iota=0}^{\infty}\frac{1}{\left(q;q\right)_{\iota}}z^{\iota}=\frac{1}{\left(z;q\right)_{\infty}}=\;_{1}\Phi_{0}\left(0;-;q,-z\right)\;\;\; \left(\left|z\right|<1\right).
\end{equation}
\end{definition}
\begin{definition}
The $q$-analogue of the Riemann–Liouville fractional integral operator for a function $f(x)$ is defined as (see, \cite{am,as,gar})
\begin{equation}  \label{1.19}
\begin{split}
	I_{q}^{\alpha}f(x)&=\frac{1}{\Gamma_{q}\left(\alpha\right)}\int_{0}^{x}\left(x-tq\right)_{\alpha-1}\,f(t)\,\mathrm{d}_{q}t\\
	&=\frac{x^{\alpha-1}}{\Gamma_{q}\left(\alpha\right)}\int_{0}^{x}\left(\frac{tq}{x};q\right)_{\alpha-1}f(t)\,\mathrm{d}_{q}t\;\;\;\; \left(\Re(\alpha)>0\right),
\end{split}
	\end{equation}
and $q$-analogue of the Riemann-Liouville fractional derivative is defined as
\begin{equation} \label{1.20}
D_{x,q}^{\alpha}f(x) = D_{q}^{n}I_{q}^{n-\alpha}f(x)\;\;\; \left(n-1<\Re(\alpha)<n,\,n\in \mathbb{N}\right)
\end{equation}
\end{definition}
\begin{definition}
The left-hand and right-hand Kober fractional $q$-integral operators are defined as (see \cite{ag})
\begin{equation} \label{1.21}
	\begin{split}
	I_{q}^{\nu,\mu}f(x)&=\frac{x^{-\nu-\mu}}{\Gamma_{q}(\mu)}\int_{0}^{x}\left(x-tq\right)_{\mu-1}t^{\nu}f(t)\,\mathrm{d}_{q}t\\
	&=\frac{x^{-\nu-1}}{\Gamma_{q}(\mu)}\int_{0}^{x}\left(\frac{tq}{x}\right)_{\mu-1}t^{\nu}f(t)\,\mathrm{d}_{q}t\;\;\; \left(\Re(\mu)>0\right),
	\end{split}
	\end{equation}
and
\begin{equation} \label{1.22}
	\begin{split}
	K_{q}^{\nu,\mu}f(x)&=\frac{x^{\nu}q^{-\nu}}{\Gamma_{q}(\mu)}\int_{x}^{\infty}\left(t-x\right)_{\mu-1}t^{-\nu-\mu}f\left(tq^{1-\mu}\right)\,\mathrm{d}_{q}t\\
	&=\frac{x^{\nu}q^{-\nu}}{\Gamma_{q}(\mu)}\int_{x}^{\infty}\left(\frac{x}{t}\right)_{\mu-1}t^{-\nu-1}f(tq^{1-\mu})\,\mathrm{d}_{q}t\;\;\; \left(\Re(\mu)>0\right),
	\end{split}
	\end{equation}
where $\mu$ is an arbitrary integration order satisfying $\Re(\mu) > 0$, and $\nu$ is a real or complex number.	
\end{definition}
\begin{definition}
The left-hand and right-hand sided Kober fractional $q$-derivative operators are given by (see \cite{gc})
\begin{equation} \label{1.23}
D_{q}^{\nu,\mu}f(x)=\prod_{i=1}^{m}\left(\left[\nu+i\right]_{q}+xq^{\nu+i}D_{q}\right)I_{q}^{\nu+\mu,m-\mu}f(x),
	\end{equation}
and
\begin{equation} \label{1.24}
\mathbb{D}_{q}^{\nu,\mu}f(x)=\prod_{i=0}^{m-1}\left(\left[\nu+i\right]_{q}-xD_{q}\right)K_{q}^{\nu+\mu,m-\mu}f(x),
\end{equation}
where $m=[\Re(\mu)]+1,\;m \in \mathbb{N}$.
\end{definition}
\begin{definition}
The basic $q$-analogue of the $q$-Weyl fractional derivative operator for an arbitrary order $\alpha$, is expressed as (see \cite{fbb})
\begin{equation} \label{1.25}
D_{\infty,q}^{\alpha}f(x)=\frac{q^{-\frac{\alpha(\alpha+1)}{2}}}{\Gamma_{q}(\alpha)}\int_{x}^{\infty}(t-x)_{-\alpha-1}f(tq^{1+\alpha})\,\mathrm{d}_{q}t\;\; \left(\Re(\alpha)>0\right).
	\end{equation}

Al-Salam \cite{al} derived the  basic $q$-analogue of the $q$-Weyl fractional integral operator for an arbitrary order $\alpha$, given by
\begin{equation} \label{1.26}
	\begin{split}
&W_{q}^{\alpha}f(x)=\frac{q^{-\frac{\alpha(\alpha-1)}{2}}}{\Gamma_{q}(\alpha)}\int_{x}^{\infty}(t-x)_{\alpha-1}f\left(tq^{1-\alpha}\right)\, \mathrm{d}_{q}t\\
&	=\frac{q^{-\frac{\alpha(\alpha-1)}{2}}}{\Gamma_{q}(\alpha)}\int_{x}^{\infty}\bigg{(}\frac{x}{t}\bigg{)}_{\alpha-1}t^{\alpha-1}f\left(tq^{1-\alpha}\right)\,\mathrm{d}_{q}t\;\; \left(\Re(\alpha)>0\right).
	\end{split}
	\end{equation}
	
Further, Al-Salam \cite{al} also introduced the $q$-extension of the generalized $q$-Weyl fractional integral operator, defined as
	\begin{equation} \label{1.27}	W_{q}^{\alpha,\beta}f(x)=\frac{q^{-\beta}x^{\beta}}{\Gamma_{q}(\alpha)}\int_{x}^{\infty}(t-x)_{\alpha-1}t^{-\alpha-\beta}f\left(tq^{1-\alpha}\right)\,\mathrm{d}_{q}t\;\; \left(\Re(\alpha)>0\right),
	\end{equation}
and $q$-analogue of $q$-Weyl fractional derivative \cite{ypk}, is defined as
\begin{equation} \label{1.28}
	\,_{-\infty}D_{x,q}^{\alpha}f(x)=(-1)^{n}D_{q}^{n}\,W_{q}^{n-\alpha}f(x)\;\; \left(n-1<\Re(\alpha)<n;\, n\in \mathbb{N}\right).
	\end{equation}
\end{definition}
\begin{definition}
Vyas et al. \cite{vap} defined the $q$-analogues of the $q$-Sumudu transform, given as
\begin{equation} \label{1.29}
	\begin{split}
	F_{s,q}(u)=S_{q}\{f(t);s\}=\frac{1}{(1-q)s}\int_{0}^{s}f(t)E_{q}\left(\frac{qt}{s}\right)\, \mathrm{d}_{q}t\\
	=\frac{1}{1-q}\int_{0}^{1}f(st)E_{q}(qt)\,\mathrm{d}_{q}t\\
	S_{q}\{t^{\alpha-1};s\}=s^{\alpha-1}\left(1-q\right)^{\alpha-1}\Gamma_{q}(\alpha)\;\;\; \left(s\in (-\tau_{1},\tau_{2})\right),
	\end{split}
	\end{equation}
for the set of functions
\begin{equation} \label{1.30}
	A=\left\{f(t)\exists M,\tau_{1}, \tau_{2}>0,\,\left|f(t)\right|<ME_{q}\left(\frac{|t|}{\tau_{j}}\right),\;t\in (-1)^{j}\times [0,\infty)\right\}.
	\end{equation}
\end{definition}
\begin{definition}
The Al-Omari \cite{als} and Nahid et al. \cite{nsc} defined the $q$-analogue of the $q$-Natural transform of the first kind over the set $A_{1}$ of the function $f$ for all $t\geq0$, given by
\begin{eqnarray} \label{1.31}
	\begin{split}
	N_{q}\left\{f(t)\right\}(u;s)&=\frac{1}{(1-q)u}\int_{0}^{\frac{u}{v}}f(t)E_{q}\left(-\frac{qst}{u}\right)\,\mathrm{d}_{q}t\;\;\; \left(s,u>0\right)\\
	&=\frac{1}{1-q}\int_{0}^{\infty}f(ut)e_{q}\left(-\frac{st}{u}\right)\,\mathrm{d}_{q}t\;\;\; \left(s,u>0\right),
	\end{split}
	\end{eqnarray}
where $t$ and $u$ are two variables, with $s$ representing the frequency variable, assuming that the function $f(t)$ is defined within a set
\begin{eqnarray} \label{1.32}
	D=\left\{f(t):\exists \tau_{1}, \tau_{2}>0, \left|f(t)\right|<M e^{\frac{t}{\tau_{j}}};\,t\in(-1)^{j}\times [0,\infty),j=1,2\right\},
		\end{eqnarray}
with $M$ as a finite constant, while $\tau_1$ and $\tau_2$ can be either finite or infinite.
\end{definition}
In this paper, the aim is to evaluate and analyze certain recurrence relations involving $q$-derivatives, $q$-Laplace transforms, $q$-Sumudu transforms, and $q$-Natural transforms, based on the series representation discussed above. In the final section, we explore the properties of the $\,{r}R{s,q}(\alpha, \beta, z)$ function using fractional $q$-calculus. We derive its $q$-fractional integrals and $q$-derivatives, and also present some special cases.

\section{Definition and convergence for $\;_{r}R_{s,q}$ function} \label{sec 2}
In this section, we discuss and evaluate the properties convergence of $\;_{r}R_{s,q}$ function.
\begin{definition}
For $r$ and $s$ are finite positive integers, and $0<|q|<1$, $q\in\mathbb{C}$, then we define the basic $\;_{r}R_{s,q}$ function by
\begin{equation} \label{2.1}
	\begin{split}
	&\;_{r}R_{s,q}\left(a_{1},a_{2},\ldots,a_{r};b_{1},b_{2},\ldots,b_{s};\alpha,\beta;z\right)\\
&	=\sum_{\ell=0}^{\infty}\frac{z^{\ell}(-1)^{\left(1+s-r\right)\ell}q^{\left(1+s-r\right)\binom{\ell}{2}}}{\left(q;q\right)_{\ell}\Gamma_{q}(\ell \alpha+\beta)}\frac{\prod_{i=1}^{r}\left(a_{i};q\right)_{\ell}}{\prod_{j=1}^{s}(b_{j};q)_{\ell}}=\sum_{\ell=0}^{\infty}U_{\ell}.
	\end{split}
	\end{equation}
\end{definition}
The expression for $U_{\ell}$ is given by
$$U_{\ell}=\frac{z^{\ell}(-1)^{((1+s-r)\ell)}q^{(1+s-r)\binom{\ell}{2}}}{(q;q)_{\ell}\Gamma_{q}(\ell \alpha+\beta)}\frac{\prod_{i=1}^{r}(a_{i};q)_{\ell}}{\prod_{j=1}^{s}(b_{j};q)_{\ell}}$$,
where it is assumed that the parameters $b_{}1$, $b_{2}$, $\ldots$, $b_{s}$ are chosen such that the denominator factors in the terms of series are never zero. Specifically, none of the
$b_{j}$'s take the values $q^{-\ell}$ for any integer $\ell=0,1,2,\ldots$.

In this paper, we consistently assume
\begin{equation}  \label{2.2}
a_{i}q^{\ell}\neq1 \; \text{and}\;  b_{j}q^{\ell}\neq1\; \left(i=1,2,\ldots,r;\,j=1,2,\ldots,s;\,\ell=0,1,2,\ldots\right),
\end{equation}
so that the factors $(a_{i};q)_{\ell}$ and $(b_{j};q)_{\ell}$ in the terms of the series are never zero.

Note that if we denote the term containing $z^{\ell}$ in \eqref{2.1} by $u_{\ell}$, then
\begin{equation} \label{2.3}
\begin{split}
\left|\frac{U_{\ell+1}}{U_{\ell}}\right|=&\left|\frac{\frac{z^{\ell+1}(-1)^{(1+s-r)(\ell+1)}q^{(1+s-r)\binom{\ell+1}{2}}}{(q;q)_{\ell+1}\Gamma_{q}(\ell \alpha+\alpha+\beta)}\frac{\prod_{i=1}^{r}(a_{i};q)_{\ell+1}}{\prod_{j=1}^{s}(b_{j};q)_{\ell+1}}}{\frac{z^{\ell}(-1)^{(1+s-r)(\ell)}q^{(1+s-r)\binom{\ell}{2}}}{(q;q)_{\ell}\Gamma_{q}(\ell \alpha+\beta)}\frac{\prod_{i=1}^{r}(a_{i};q)_{\ell}}{\prod_{j=1}^{s}(b_{j};q)_{\ell}}}\right|\\
=&\left|\frac{z(-1)^{(1+s-r)}q^{(1+s-r)(\ell)}\Gamma_{q}(\ell \alpha+\beta)}{(1-q^{\ell+1})\Gamma_{q}(\ell \alpha+\alpha+\beta)}\frac{\prod_{i=1}^{r}(1-a_{i}q^{\ell})}{\prod_{j=1}^{s}(1-b_{j}q^{\ell})}\right|\\
=&\left|\frac{\Gamma_{q}(\ell \alpha+\beta)}{(1-q^{-\ell-1})\Gamma_{q}(\ell \alpha+\alpha+\beta)}\frac{\prod_{i=1}^{r}(a_{i}-q^{-\ell})}{\prod_{j=1}^{s}(b_{j}-q^{-\ell})}\frac{z}{q}\right|.
\end{split}
\end{equation}
Thus, by the ratio test we can say that
\begin{enumerate}
\item If $0<|q|<1$, $q\in\mathbb{C}$, $s+1\geq r$, the $\;_{r}R_{s,q}$ series converges absolutely for all $z \in \mathbb{C}$.
\item If $0<|q|<1$, $q\in\mathbb{C}$, $s+2<r$, the $\;_{r}R_{s,q}$ series diverges for $z\neq 0$.
\item If $0<|q|<1$, $q\in\mathbb{C}$, $s+2=r$, the $\;_{r}R_{s,q}$ series converges for $|\frac{z}{q}|<1$ and diverges for $|\frac{z}{q}|>1$.
\item If $|q|>1$, $q\in\mathbb{C}$, $s+2=r$, the  $\;_{r}R_{s,q}$ series converges absolutely for $|z|<\frac{|b_{1}b_{2}\ldots b_{s}q|}{|a_{1}a_{2} \ldots a_{r}|}$ and diverges for $|z|>\frac{|b_{1}b_{2}\ldots b_{s}q|}{|a_{1}a_{2}\ldots a_{r}|}$.
\end{enumerate}
\begin{corollary}
The limit relation between $\;_{r}R_{s}$ series and  $\;_{r}R_{s,q}$ series is given by
\begin{equation} \label{2.4}
	\begin{split}
	&\lim_{q\Rightarrow 1}\;_{r}R_{s,q}\left(q^{a_{1}},q^{a_{2}},\ldots,q^{a_{r}};q^{b_{1}},q^{b_{2}},\ldots,q^{b_{s}};\alpha,\beta;\frac{(-1)^{r-s-1}z}{(1-q)^{r-s-1}}\right)\\
&=\;_{r}R_{s}(a_{1},a_{2},\ldots,a_{r};b_{1},b_{2},\ldots,b_{s};\alpha,\beta;z).
	\end{split}
	\end{equation}
\end{corollary}
\begin{proof}
By the definition of $q$-number and $q$-shifted factorial and taking the limit of both sides of this equation when $q$ approaches 1, we get
\begin{equation} \label{2.5}
\lim_{q\Rightarrow 1}\frac{(q^{a_{i}};q)_{\ell}}{\left(1-q\right)^{\ell}}=\left(a_{i}\right)_{\ell}.
\end{equation}
Alternatively, using the definition of the  $\;_{r}R_{s,q}$ series, we obtain
\begin{equation*}
	\begin{split}
	&\lim_{q\Rightarrow 1}\;_{r}R_{s,q}\left(q^{a_{1}},q^{a_{2}},\ldots,q^{a_{r}};q^{b_{1}},q^{b_{2}},\ldots,q^{b_{s}};\alpha,\beta;\frac{(-1)^{r-s-1}z}{(1-q)^{r-s-1}}\right)\\
	&=\;_{r}R_{s}\left(a_{1},a_{2},\ldots,a_{r};\, b_{1},b_{2},\ldots,b_{s};\,\alpha,\beta;z\right).
	\end{split}
	\end{equation*}
This yields the desired result.
\end{proof}

\section{Recurrence relations}   \label{sec 3}
In this section, we establish the recurrence relations and differential properties for the $\;_{r}R_{s,q}$ function by using fractional $q$-calculus.
\begin{theorem}
The basic hypergeometric function is satisfied by the identities
\begin{equation} \label{3.1}
\left(1-a_{i}\right)\;_{r}R_{s,q}(a_{i}q)=\;_{r}R_{s,q}-a_{i}\;_{r}R_{s,q}(qz)=0\,\; \left(i=1,2,\ldots,r\right),
	\end{equation}
	\begin{equation}   \label{3.2}
	\;_{r}R_{s,q}(a_{1}q)=\;_{r}R_{s,q}+(-q)^{1+s-r}a_{1}z\frac{\prod_{i=2}^{r}(1-a_{i})}{\prod_{j=1}^{s}(1-b_{j})}\;_{r}R_{s,q}\left(a_{i}q;b_{j}q;\alpha,\alpha+\beta,zq^{1+s-r}\right),
	\end{equation}
	\begin{equation} \label{3.3}
\left(1-b_{j}q^{-1}\right)\;_{r}R_{s,q}\left(b_{j}q^{-1}\right)=\;_{r}R_{s,q}-b_{j}q^{-1}\;_{r}R_{s,q}\left(qz\right)=0,
	\end{equation}
and
\begin{equation} \label{3.4}
	\;_{r}R_{s,q}(b_{1}q^{-1})=\,_{r}R_{s,q}+\frac{(-q)^{1+s-r}b_{1}q^{-1}}{1-b_{1}q^{-1}}z\frac{\prod_{i=1}^{r}(1-a_{i})}{\prod_{j=1}^{s}(1-b_{j})}\;_{r}R_{s,q}\left(a_{i}q;b_{j}q;\alpha,\alpha+\beta,zq^{1+s-r}\right),
	\end{equation}
where, $i=1,2,\ldots,r;\,j=1,2,\ldots,s$.
\end{theorem}
\begin{proof}
By using the relation
\begin{equation} \label{3.5}
\left(1-a_{i}\right)\left(a_{i}q;q\right)_{\ell}=\left(1-a_{i}q^{\ell}\right)\left(a_{i};q\right)_{\ell}=\left(a_{i};q\right)_{\ell+1}.
\end{equation}
To prove \eqref{3.1}, we proceed by applying \eqref{3.5}, which gives
\begin{equation*}
	\begin{split}
	\;_{r}R_{s,q}(a_{1}q)-\;_{r}R_{s,q}=&\sum_{\ell=0}^{\infty}\frac{z^{\ell}(-1)^{(1+s-r)\ell}q^{(1+s-r)\binom{\ell}{2}}}{(q;q)_{\ell}\Gamma_{q}(\ell \alpha+\beta)}\frac{\prod_{i=2}^{r}(a_{i};q)_{\ell}}{\prod_{j=1}^{s}\left(b_{j};q\right)_{\ell}}\left[\left(a_{1}q;q\right)_{\ell}-\left(a_{1};q\right)_{\ell}\right]\\
	=&\sum_{\ell=0}^{\infty}\frac{z^{\ell}(-1)^{(1+s-r)\ell}q^{(1+s-r)\binom{\ell}{2}}}{\left(q;q\right)_{\ell}\Gamma_{q}\left(\ell \alpha+\beta\right)}\frac{\prod_{i=2}^{r}(a_{i};q)_{\ell}}{\prod_{j=1}^{s}(b_{j};q)_{\ell}}\left(a_{1};q\right)_{\ell}\left[\frac{1-a_{1}q^{\ell}}{1-a_{1}}-1\right]\\
	=&\frac{a_{1}}{1-a_{1}}\sum_{\ell=0}^{\infty}\frac{\left(1-q^{\ell}\right)z^{\ell}(-1)^{\left(1+s-r\right)\ell}q^{\left(1+s-r\right)\binom{\ell}{2}}}{\left(q;q\right)_{\ell}\Gamma_{q}\left(\ell \alpha+\beta\right)}\frac{\prod_{i=2}^{r}\left(a_{i};q\right)_{\ell}}{\prod_{j=1}^{s}\left(b_{j};q\right)_{\ell}}\left(a_{1};q\right)_{\ell}\\
	=&\frac{a_{1}}{1-a_{1}}\left[\;_{r}R_{s,q}-\;_{r}R_{s,q}(qz)\right].
	\end{split}
	\end{equation*}
Making use of the identities given in \eqref{3.5}, we obtain
\begin{equation*}
	\begin{split}
	&\;_{r}R_{s,q}(a_{1}q)-\;_{r}R_{s,q}\\
	=&\sum_{\ell=0}^{\infty}\frac{z^{\ell}(-1)^{(1+s-r)\ell}q^{(1+s-r)\binom{\ell}{2}}}{(q;q)_{\ell}\Gamma_{q}\left(\ell \alpha+\beta\right)}\frac{\prod_{i=2}^{r}\left(a_{i};q\right)_{\ell}}{\prod_{j=1}^{s}\left(b_{j};q\right)_{\ell}}\left[\left(a_{1}q;q\right)_{\ell}-\frac{1-a_{1}}{1-a_{1}q^{\ell}}\left(a_{1}q;q\right)_{\ell}\right]\\
	=&\sum_{\ell=0}^{\infty}\frac{z^{\ell}(-1)^{\left(1+s-r\right)\ell}q^{\left(1+s-r\right)\binom{\ell}{2}}}{\left(q;q\right)_{\ell}\Gamma_{q}\left(\ell \alpha+\beta\right)}\frac{\prod_{i=2}^{r}\left(a_{i};q\right)_{\ell}}{\prod_{j=1}^{s}\left(b_{j};q\right)_{\ell}}\left(a_{1}q;q\right)_{\ell}\left[\frac{a_{1}(1-q^{\ell})}{1-a_{1}q^{\ell}}\right]\\
	=&\sum_{\ell=0}^{\infty}\frac{z^{\ell+1}(-1)^{\left(1+s-r\right)\left(\ell+1\right)}q^{\left(1+s-r\right)\binom{\ell+1}{2}}}{\left(q;q\right)_{\ell}\Gamma_{q}\left(\ell \alpha+\alpha+\beta\right)}\frac{\prod_{i=2}^{r}\left(1-a_{i}\right)\left(a_{i}q;q\right)_{\ell}}{\prod_{j=1}^{s}\left(1-b_{j}\right)\left(b_{j}q;q\right)_{\ell}}a_{1}\left(a_{1}q;q\right)_{\ell}\\
	=&(-q)^{1+s-r}a_{1}z\frac{\prod_{i=2}^{r}(1-a_{i})}{\prod_{j=1}^{s}(1-b_{j})}\;_{r}R_{s,q}\left(a_{1}q,a_{2}q,...,a_{r}q;\,b_{1}q,b_{2}q,\ldots,b_{s}q;\,\alpha,\alpha+\beta,\,zq^{1+s-r}\right).
	\end{split}
	\end{equation*}
By replacing $b_{1}$ with $b_{1}q^{-1}$ in \eqref{2.1}, we obtain
\begin{equation*}
	\begin{split}
	\;_{r}R_{s,q}(b_{1}q^{-1})-\,_{r}R_{s,q}=&\sum_{\ell=0}^{\infty}\frac{z^{\ell}(-1)^{(1+s-r)\ell}q^{\left(1+s-r\right)\binom{\ell}{2}}}{\left(q;q\right)_{\ell}\,\Gamma_{q}\left(\ell \alpha+\beta\right)}\frac{\prod_{i=1}^{r}(a_{i};q)_{\ell}}{\prod_{j=2}^{s}\left(b_{j};q\right)_{\ell}}\left[\frac{1}{\left(b_{1}q^{-1};q\right)_{\ell}}-\frac{1}{\left(b_{1};q\right)_{\ell}}\right].
	\end{split}
	\end{equation*}
By using the following relation:
\begin{equation} \label{3.6}	\frac{1}{\left(b_{1};q\right)_{\ell-1}}=\frac{\left(1-b_{1}q^{\ell-1}\right)}{\left(b_{1};q\right)_{\ell}}=\frac{\left(1-b_{1}q^{-1}\right)}{\left(b_{1}q^{-1};q\right)_{\ell}},
	\end{equation}
	we have
	\begin{equation*}
	\begin{split}
&	\;_{r}R_{s,q}(b_{1}q^{-1})-\;_{r}R_{s,q}\\
&=\sum_{\ell=0}^{\infty}\frac{z^{\ell}(-1)^{\left(1+s-r\right)\ell}q^{\left(1+s-r\right)\binom{\ell}{2}}}{\left(q;q\right)_{\ell}\Gamma_{q}\left(\ell \alpha+\beta\right)}\frac{\prod_{i=1}^{r}\left(a_{i};q\right)_{\ell}}{\prod_{j=2}^{s}\left(b_{j};q\right)_{\ell}}\left[\frac{\left(b_{1};q\right)_{\ell}-\left(b_{1}q^{-1};q\right)_{\ell}}{\left(b_{1};q\right)_{\ell}\left(b_{1}q^{-1};q\right)_{\ell}}\right]\\
&	=\sum_{\ell=0}^{\infty}\frac{z^{\ell}(-1)^{\left(1+s-r\right)\ell}q^{\left(1+s-r\right)\binom{\ell}{2}}}{\left(q;q\right)_{\ell-1}\Gamma_{q}\left(\ell \alpha+\beta\right)}\frac{\prod_{i=1}^{r}\left(a_{i};q\right)_{\ell}}{\prod_{j=2}^{s}\left(b_{j};q\right)_{\ell}}\left[\frac{b_{1}q^{-1}}{\left(1-b_{1}q^{-1}\right)\left(b_{1};q\right)_{\ell}}\right]\\
&	=\frac{b_{1}\prod_{i=1}^{r}(1-a_{i})}{\left(1-b_{1}\right)\left(q-b_{1}\right)\prod_{j=2}^{s}\left(1-b_{j}\right)}\sum_{\ell=0}^{\infty}\frac{z^{\ell}(-1)^{\left(1+s-r\right)\ell}q^{\left(1+s-r\right)\binom{\ell}{2}}}{\left(q;q\right)_{\ell-1}\Gamma_{q}\left(\ell \alpha+\beta\right)}\frac{\prod_{i=1}^{r}\left(a_{i}q;q\right)_{\ell-1}}{\left(b_{1}q;q\right)_{\ell-1}\prod_{j=2}^{s}\left(b_{j}q;q\right)_{\ell-1}}\\
&	=\frac{b_{1}z(-q)^{1+s-r}\prod_{i=1}^{r}\left(1-a_{i}\right)}{\left(1-b_{1}\right)\left(q-b_{1}\right)\prod_{j=2}^{s}\left(1-b_{j}\right)}\;_{r}R_{s,q}\left(a_{i}q,a_{2}q,\ldots,a_{r}q;b_{1}q,b_{2}q,\ldots,b_{s}q;\alpha,\alpha+\beta,zq^{1+s-r}\right),
	\end{split}
	\end{equation*}
We prove in a similar manner that of \eqref{3.3} and \eqref{3.4}.
\end{proof}
\begin{theorem}
Each of the following relations for $\;_{r}R_{s,q}$ hold true:
\begin{equation} \label{3.7}
	\begin{split}
	&\;_{r}R_{s,q}(a_{1}q,b_{1}q)-\,_{r}R_{s,q}=\frac{a_{1}z(-q)^{1+s-r}\left(1-\frac{b_{1}}{a_{1}}\right)}{1-b_{1}q}\\
	&\times\;\frac{\prod_{i=2}^{r}(1-a_{i})}{\prod_{j=1}^{s}\left(1-b_{j}\right)}\;_{r}R_{s,q}\left(a_{1}q,a_{2}q,\ldots;a_{r}q,\;b_{1}q^{2},b_{2}q,\ldots,b_{s}q,\;\alpha,\alpha+\beta;\;zq^{1+s-r}\right),
	\end{split}
	\end{equation}
and
\begin{equation} \label{3.8}
	\begin{split}
	&\;_{r}R_{s,q}\left(a_{1}q,a_{2}q^{-1}\right)-\;_{r}R_{s,q}=z(-q)^{1+s-r}\left(a_{1}-a_{2}q^{-1}\right)\\
	&\times\,\frac{\prod_{i=3}^{r}\left(1-a_{i}\right)}{\prod_{j=1}^{s}\left(1-b_{j}\right)}\;_{r}R_{s,q}\left(a_{1}q,a_{2};a_{3}q,\ldots;a_{r}q,\,b_{1}q^{2},b_{2}q,\ldots,b_{s}q,\,\alpha,\alpha+\beta;\,zq^{1+s-r}\right).
	\end{split}
	\end{equation}
\end{theorem}
\begin{proof}
	In equation \eqref{2.1}, replacing $a$ and $b$ by $aq$ and $bq$, respectively, and by using relation \eqref{3.5}, we arrive at
\begin{equation*}
	\begin{split}
	&\;_{r}R_{s,q}(a_{1}q,b_{1}q)-\;_{r}R_{s,q}\\
&=\sum_{\ell=0}^{\infty}\frac{z^{\ell}(-1)^{(1+s-r)\ell}q^{(1+s-r)\binom{\ell}{2}}}{(q;q)_{\ell}\Gamma_{q}\left(\ell \alpha+\beta\right)}\frac{\prod_{i=2}^{r}\left(a_{i};q\right)_{\ell}}{\prod_{j=2}^{s}\left(b_{j};q\right)_{\ell}}\left[\frac{\left(a_{1}q;q\right)_{\ell}\left(b_{1};q\right)_{\ell}-\left(b_{1}q;q\right)_{\ell}\left(a_{1};q\right)_{\ell}}{\left(b_{1}q;q\right)_{\ell}\left(b_{1};q\right)_{\ell}}\right]\\
	&=\sum_{\ell=0}^{\infty}\frac{z^{\ell}(-1)^{\left(1+s-r\right)\ell}q^{(1+s-r)\binom{\ell}{2}}\left(a_{1}q;q\right)_{\ell-1}\left(b_{1}q;q\right)_{\ell-1}}{\left(q;q\right)_{\ell}\Gamma_{q}\left(\ell \alpha+\beta\right)}\\
&\times\;	\frac{\prod_{i=2}^{r}(a_{i};q)_{\ell}}{\prod_{j=2}^{s}\left(b_{j};q\right)_{\ell}}\left[\frac{\left(a_{1}-b_{1}\right)\left(1-q^{\ell}\right)}{\left(1-b_{1}q\right)\left(b_{1}q;q\right)_{\ell-1}\left(1-b_{1}\right)\left(b_{1}q;q\right)_{\ell-1}}\right]\\
&=\frac{a_{1}\left(1-\frac{b_{1}}{a_{1}}\right)}{\left(1-b_{1}\right)\left(1-b_{1}q\right)}\sum_{\ell=0}^{\infty}\frac{z^{\ell}(-1)^{\left(1+s-r\right)\ell}q^{\left(1+s-r\right)\binom{\ell}{2}}\left(a_{1}q;q\right)_{\ell-1}}{\left(q;q\right)_{\ell-1}\left(b_{1}q^{2};q\right)_{\ell-1}\Gamma_{q}\left(\ell \alpha+\beta\right)}\frac{\prod_{i=2}^{r}\left(1-a_{i}\right)\left(a_{i}q;q\right)_{\ell-1}}{\prod_{j=2}^{s}\left(1-b_{j}\right)\left(b_{j}q;q\right)_{\ell-1}}.
\end{split}
\end{equation*}
If write $\ell+1$ in place of $\ell$, it follows that
\begin{equation*}
	\begin{split}
	&\;_{r}R_{s,q}\left(a_{1}q,b_{1}q\right)-\,_{r}R_{s,q}=\frac{a_{1}\left(1-\frac{b_{1}}{a_{1}}\right)}{\left(1-b_{1}\right)\left(1-b_{1}q\right)}\sum_{\ell=-1}^{\infty}\frac{z^{\ell+1}(-1)^{\left(1+s-r\right)\left(\ell+1\right)}q^{\left(1+s-r\right)\binom{\ell+1}{2}}\left(a_{1}q;q\right)_{\ell}}{\left(q;q\right)_{\ell}\left(b_{1}q^{2};q\right)_{\ell}\Gamma_{q}\left(\ell \alpha+\alpha+\beta\right)}\\
& \times\; \frac{\prod_{i=2}^{r}\left(1-a_{i}\right)\left(a_{i}q;q\right)_{\ell}}{\prod_{j=2}^{s}\left(1-b_{j}\right)\left(b_{j}q;q\right)_{\ell}} =\frac{(-q)^{1+s-r}za_{1}\left(1-\frac{b_{1}}{a_{1}}\right)}{\left(1-b_{1}\right)\left(1-b_{1}q\right)}\frac{\prod_{i=2}^{r}\left(1-a_{i}\right)}{\prod_{j=2}^{s}\left(1-b_{j}\right)}\\
& \times\; \;_{r}R_{s,q}\left(a_{1}q,a_{2}q,\ldots,a_{r}q,\,b_{1}q^{2},b_{2}q,\ldots,b_{s}q;\,\alpha,\alpha+\beta,\,zq^{1+s-r}\right).
	\end{split}
	\end{equation*}
From equations \eqref{2.1}, \eqref{3.5} and \eqref{3.6}, we can write
\begin{equation*}
	\begin{split}
	&\,_{r}R_{s,q}(a_{1}q,a_{2}q^{-1})-\,_{r}R_{s,q}=\sum_{\ell=0}^{\infty}\frac{z^{\ell}(-1)^{(1+s-r)\ell}q^{(1+s-r)\binom{\ell}{2}}}{(q;q)_{\ell}\Gamma_{q}(\ell \alpha+\beta)}\frac{(a_{1}q;q)_{\ell-1}(a_{2};q)_{\ell-1}\prod_{i=3}^{r}(a_{i};q)_{\ell}}{\prod_{j=1}^{s}(b_{j};q)_{\ell}}\\
	&\times\,\left[\left(1-a_{1}q^{\ell}\right)\left(1-a_{2}q^{-1}\right)-\left(1-a_{1}\right)\left(1-a_{2}q^{\ell-1}\right)\right]\\
	&=\left(a_{1}-a_{2}q^{-1}\right)\sum_{\ell=0}^{\infty}\frac{z^{\ell}(-1)^{(1+s-r)\ell}q^{(1+s-r)\binom{\ell}{2}}}{(q;q)_{\ell-1}\Gamma_{q}(\ell \alpha+\beta)}\;\frac{(a_{1}q;q)_{\ell-1}(a_{2};q)_{\ell-1}\prod_{i=3}^{r}(a_{i};q)_{\ell}}{\prod_{j=1}^{s}(b_{j};q)_{\ell}}\\
	&=\left(a_{1}-a_{2}q^{-1}\right)\sum_{\ell=0}^{\infty}\frac{z^{\ell+1}(-1)^{(1+s-r)(\ell+1)}q^{(1+s-r)\binom{\ell+1}{2}}}{(q;q)_{\ell}\Gamma_{q}(\ell \alpha+\alpha+\beta)}\,\frac{(a_{1}q;q)_{\ell}(a_{2};q)_{\ell}\prod_{i=3}^{r}(1-a_{i})(a_{i}q;q)_{\ell}}{\prod_{j=1}^{s}(1-b_{j})(b_{j}q;q)_{\ell}}\\
	&=\frac{(-q)^{1+s-r}z(a_{1}-a_{2}q^{-1})\prod_{i=3}^{r}(1-a_{i})}{\prod_{j=1}^{s}(1-b_{j})}\\
&\times\;	\;_{r}R_{s,q}\left(a_{1}q,a_{2},\ldots,a_{r}q,\,b_{1}q,b_{2}q,\ldots,b_{s}q;\,\alpha,\alpha+\beta,\,zq^{1+s-r}\right).
	\end{split}
\end{equation*}
\end{proof}
\begin{theorem}
The $\;_{r}R_{s,q}$ function have the recurrence relations
\begin{equation} \label{3.9}
a_{k}\left(1-a_{i}\right)\,_{r}R_{s,q}\left(a_{i}q\right)-a_{i}\left(1-a_{k}\right)\;_{r}R_{s,q}\left(a_{k}q\right)=\left(a_{k}-a_{i}\right)\,_{r}R_{s,q}\; \left(i\neq k, i,k=1,2,\ldots,r\right),
	\end{equation}
and
\begin{equation} \label{3.10}
	\begin{split}
&	b_{k}q^{-1}\left(1-b_{j}q^{-1}\right)\,_{r}R_{s,q}\left(b_{k}q^{-1}\right)-b_{j}q^{-1}\left(1-b_{k}q^{-1}\right)\,_{r}R_{s,q}\left(b_{j}q^{-1}\right)\\
&	=\left(b_{k}-b_{j}\right)q^{-1}\,_{r}R_{s,q}(b_{j}q^{-1},b_{k}q^{-1})\;\;\; \left(i=1,2,\ldots,r;\,k=1,2,\ldots,s,\,i\neq k\right).
	\end{split}
	\end{equation}
\end{theorem}
\begin{proof}
For $i\neq k$, $i,k=1,2,\ldots,r$, we get the following relation:
\begin{equation} \label{3.11}
	\begin{split}
&	\left(1-a_{i}\right)\left(a_{i}q;q\right)_{\ell}=\left(1-a_{i}q^{\ell}\right)\left(a_{i};q\right)_{\ell},\\
&	\left(1-a_{k}\right)\left(a_{k}q;q\right)_{\ell}=\left(1-a_{k}q^{\ell}\right)\left(a_{k};q\right)_{\ell}.
	\end{split}
	\end{equation}
Again, by making use of the equation \eqref{3.5} and the above relations as stated in \eqref{3.11}, yields the result
\begin{equation*}
	\begin{split}
	&a_{2}(1-a_{1})\;_{r}R_{s,q}(a_{1}q)-a_{1}\left(1-a_{2}\right)\;_{r}R_{s,q}\left(a_{2}q\right)\\
	&=a_{2}\left(1-a_{1}\right)\sum_{\ell=0}^{\infty}\frac{z^{\ell}(-1)^{\left(1+s-r\right)\ell}q^{\left(1+s-r\right)\binom{\ell}{2}}}{\left(q;q\right)_{\ell}\Gamma_{q}\left(\ell \alpha+\beta\right)}\,\frac{\left(a_{1}q;q\right)_{\ell}\prod_{i=2}^{r}\left(a_{i};q\right)_{\ell}}{\prod_{j=1}^{s}\left(b_{j};q\right)_{\ell}}\\
	&-a_{1}\left(1-a_{2}\right)\sum_{\ell=0}^{\infty}\,\frac{z^{\ell}(-1)^{\left(1+s-r\right)\ell}q^{\left(1+s-r\right)\binom{\ell}{2}}}{\left(q;q\right)_{\ell}\Gamma_{q}\left(\ell \alpha+\beta\right)}\,\frac{\left(a_{2}q;q\right)_{\ell}\prod_{i=2}^{r}\left(a_{i};q\right)_{\ell}}{\prod_{j=1}^{s}\left(b_{j};q\right)_{\ell}}.
	\end{split}
	\end{equation*}
By using the above equation, we get
\begin{equation} \label{3.12}
	\begin{split}
&	a_{k}\left(1-a_{i}\right)\left(a_{i}q;q\right)_{\ell}\left(a_{k};q\right)_{\ell}-a_{i}\left(1-a_{k}\right)\left(a_{i};q\right)_{\ell}\left(a_{k}q;q\right)_{\ell}\\
&=\left[a_{k}-a_{i}\right]\left(a_{i};q\right)_{\ell}\left(a_{k};q\right)_{\ell},
	\end{split}
	\end{equation}
by using \eqref{2.1} and \eqref{3.12}, we derive \eqref{3.10}.
\end{proof}
\begin{theorem}
The $\;_{r}R_{s,q}$ function has the relation
\begin{equation} \label{3.13}
	\begin{split}
&	z\frac{\prod_{i=1}^{r}\left(1-a_{i}\right)}{\prod_{j=1}^{s}\left(1-b_{j}\right)}\;_{r}R_{s,q}\left(a_{i}q,b_{j}q\right)\\
&=(-q)^{1+s-r}\,_{r}R_{s,q}\left(\alpha,\beta-\alpha;q,zq^{r-s-1}\right)-(-q)^{1+s-r}\,_{r}R_{s,q}\left(\alpha,\beta-\alpha;q,zq^{r-s}\right).
	\end{split}
	\end{equation}
\end{theorem}
\begin{proof}
From the relations \eqref{3.5} and \eqref{3.6}, we obtain
\begin{equation} \label{3.14}
	\begin{split}
\frac{\left(1-a_{i}\right)\left(a_{i}q;q\right)_{\ell}}{\left(1-b_{j}\right)\left(b_{j}q;q\right)_{\ell}}=\frac{\left(1-a_{i}q^{\ell}\right)\left(a_{i};q\right)_{\ell}}{\left(1-b_{j}q^{\ell}\right)\left(b_{j};q\right)_{\ell}}=\frac{\left(a_{i};q\right)_{\ell+1}}{\left(b_{j};q\right)_{\ell+1}}.
	\end{split}
	\end{equation}
By using of equation \eqref{3.14}, we have
\begin{equation*}
\begin{split}
&z\frac{\prod_{i=1}{r}\left(1-a_{i}\right)}{\prod_{j=1}{s}\left(1-b_{j}\right)}\;_{r}R_{s,q}\left(a_{i}q,b_{j}q\right)=\sum_{\ell=0}^{\infty}\frac{z^{\ell+1}(-1)^{\left(1+s-r\right)\ell}q^{\left(1+s-r\right)\binom{\ell}{2}}}{\left(q;q\right)_{\ell}\Gamma_{q}\left(\ell \alpha+\beta\right)}\frac{\prod_{i=1}^{r}\left(a_{i}q;q\right)_{\ell+1}}{\prod_{j=1}^{s}\left(b_{j}q;q\right)_{\ell+1}}\\
&=\sum_{\ell=0}^{\infty}\frac{z^{\ell}(-1)^{\left(1+s-r\right)\left(\ell-1\right)}q^{\left(1+s-r\right)\binom{\ell-1}{2}}}{\left(q;q\right)_{\ell-1}\Gamma_{q}\left(\ell \alpha+\beta-\alpha\right)}\frac{\prod_{i=1}^{r}\left(a_{i}q;q\right)_{\ell}}{\prod_{j=1}^{s}\left(b_{j}q;q\right)_{\ell}}\\
&=\sum_{\ell=0}^{\infty}\frac{z^{\ell}\left(1-q^{\ell}\right)(-1)^{\left(1+s-r\right)\left(\ell-1\right)}q^{\left(1+s-r\right)\binom{\ell-1}{2}}}{\left(q;q\right)_{\ell}\Gamma_{q}\left(\ell \alpha+\beta-\alpha\right)}\frac{\prod_{i=1}^{r}\left(a_{i}q;q\right)_{\ell}}{\prod_{j=1}^{s}\left(b_{j}q;q\right)_{\ell}}.
\end{split}
\end{equation*}
\end{proof}
\begin{corollary}
The following relation holds true:
\begin{equation} \label{3.15}
	\;_{r}R_{s,q}(qz)=\;_{r}R_{s,q}-z(-q)^{1+s-r}\frac{\prod_{i=1}^{r}\left(1-a_{i}\right)}{\prod_{j=1}^{s}\left(1-b_{j}\right)}\;_{r}R_{s,q}\left(a_{i}q;b_{j}q;\alpha,\alpha+\beta,zq^{1+s-r}\right)=0.
\end{equation}
\end{corollary}
\begin{theorem}
The following formulas hold true:
\begin{equation} \label{3.16}
\begin{split}
&\;_{r}R_{s,q}\left(a_{1},a_{2},\ldots,a_{r};b_{1},b_{2},\ldots,b_{s};\alpha,\beta;z\right)=\, _{r}R_{s,q}\left(a_{i}q;b_{1},b_{2},\ldots,b_{s};\alpha,\beta;z\right)\\
&	-\frac{(-q)^{1+s-r}a_{i}z}{\prod_{j=1}^{s}\left(1-b_{j}\right)}\;_{r}R_{s,q}\left(a_{i}q;b_{1}q,b_{2}q,\ldots,b_{s}q;\alpha,\alpha+\beta;q^{1+s-r}z\right)\;\; \left(i=1,2,\ldots,r\right),
	\end{split}
	\end{equation}
and
\begin{equation} \label{3.17}
	\begin{split}
	&\;_{r}R_{s,q}\left(a_{1},a_{2},\ldots,a_{r};b_{1},b_{2},\ldots,b_{s};\alpha,\beta;z\right)=\, _{r}R_{s,q}\left(a_{i}q;b_{1},b_{2},\ldots,b_{s};\alpha,\beta;z\right)\\
&	-\frac{(-q)^{1+s-r}a_{i}z}{\prod_{j=1}^{s}\left(1-b_{j}\right)}\;_{r}R_{s,q}\left(a_{i}q;b_{1}q,b_{2}q,\ldots,b_{s}q;\alpha,\alpha+\beta;q^{1+s-r}z\right)\\
&	+\frac{(-q)^{1+s-r}a_{i}qz}{\prod_{j=1}^{s}\left(1-b_{j}\right)}\, _{r}R_{s,q}\left(a_{i}q;b_{1}q,b_{2}q,\ldots,b_{s}q;\alpha,\alpha+\beta;q^{2+s-r}z\right)\;\;\left(i=1,2,\ldots,r\right).
	\end{split}
	\end{equation}
\end{theorem}
\begin{proof}
From \eqref{2.1} and use the relation $1-a_{i}=1-a_{i}q^{\ell}-a_{i}(1-q^{\ell})$, we can write
\begin{equation*}
\begin{split}
	&\;_{r}R_{s,q}\left(a_{1},a_{2},\ldots,a_{r};b_{1},b_{2},\ldots,b_{s};\alpha,\beta;z\right)\\
&	=\frac{1}{\Gamma_{q}(\beta)}+\sum_{\ell=1}^{\infty}\frac{z^{\ell}(-1)^{(1+s-r)(\ell)}q^{(1+s-r)\binom{\ell}{2}}}{\left(q;q\right)_{\ell}\Gamma_{q}\left(\ell \alpha+\beta\right)}\frac{\prod_{i=1}^{r}(a_{i};q)_{\ell}}{\prod_{j=1}^{s}\left(b_{j};q\right)_{\ell}}\\
&	=\frac{1}{\Gamma_{q}(\beta)}+\sum_{\ell=1}^{\infty}\frac{z^{\ell}(-1)^{(1+s-r)(\ell)}q^{(1+s-r)\binom{\ell}{2}}}{\left(q;q\right)_{\ell}\Gamma_{q}\left(\ell \alpha+\beta\right)}\frac{\prod_{i=1}^{r}\left(1-a_{i}q^{\ell}\right)\left(a_{i}q;q\right)_{\ell-1}}{\prod_{j=1}^{s}\left(1-b_{j}\right)\left(b_{j}q;q\right)_{\ell-1}}\\
&	-\sum_{\ell=1}^{\infty}\frac{z^{\ell}(-1)^{\left(1+s-r\right)\left(\ell\right)}q^{(1+s-r)\binom{\ell}{2}}}{\left(q;q\right)_{\ell}\Gamma_{q}\left(\ell \alpha+\beta\right)}\frac{\prod_{i=1}^{r}\left(a_{i}\left(1-q^{\ell}\right)\right)\left(a_{i}q;q\right)_{\ell-1}}{\prod_{j=1}^{s}\left(1-b_{j}\right)\left(b_{j}q;q\right)_{\ell-1}}\\
&	-a_{i}\sum_{\ell=1}^{\infty}\frac{z^{\ell}(-1)^{\left(1+s-r\right)\left(\ell\right)}q^{\left(1+s-r\right)\binom{\ell}{2}}}{\left(q;q\right)_{\ell}\Gamma_{q}\left(\ell \alpha+\beta\right)}\frac{\prod_{i=1}^{r}\left(1-q^{\ell}\right)\left(a_{i}q;q\right)_{\ell-1}}{\prod_{j=1}^{s}\left(1-b_{j}\right)\left(b_{j}q;q\right)_{\ell-1}}\\
&	=\sum_{\ell=0}^{\infty}\frac{z^{\ell}(-1)^{(1+s-r)\left(\ell\right)}q^{\left(1+s-r\right)\binom{\ell}{2}}}{\left(q;q\right)_{\ell}\Gamma_{q}\left(\ell \alpha+\beta\right)}\frac{\prod_{i=1}^{r}\left(a_{i}q;q\right)_{\ell}}{\prod_{j=1}^{s}\left(1-b_{j}\right)\left(b_{j}q;q\right)_{\ell-1}}\\
&	-a_{i}\sum_{\ell=0}^{\infty}\frac{z^{\ell+1}(-1)^{\left(1+s-r\right)\left(\ell+1\right)}q^{\left(1+s-r\right)\binom{\ell+1}{2}}}{\left(q;q\right)_{\ell}\Gamma_{q}\left(\ell \alpha+\beta+\alpha\right)}\frac{\prod_{i=1}^{r}\left(a_{i}q;q\right)_{\ell}}{\prod_{j=1}^{s}\left(1-b_{j}\right)\left(b_{j}q;q\right)_{\ell}}\\
&	=\,_{r}R_{s,q}\left(a_{i}q;b_{1},b_{2},\ldots,b_{s};\alpha,\beta;z\right)\\
&	-\frac{(-q)^{1+s-r}a_{i}z}{\prod_{j=1}^{s}\left(1-b_{j}\right)}\,_{r}R_{s,q}\left(a_{i}q;b_{1}q,b_{2}q,\ldots,b_{s}q;\alpha,\alpha+\beta;q^{1+s-r}z\right)\;\; \left(i=1,2,\ldots,r\right).
	\end{split}
	\end{equation*}
Similarly, the result in \eqref{3.17} can be proven.
\end{proof}

\begin{theorem}
The following $q$-derivative formulas for $\;_{r}R_{s,q}$ function holds true:
\begin{equation} \label{3.18}
\begin{split}
&\theta_{z,q}^{n}\;_{r}R_{s,q}=\left(zD_{z,q}\right)^{n}\,_{r}R_{s,q} = \frac{z^{n}(-q)^{\left(1+s-r\right)n}\prod_{i=1}^{r}\left(a_{i};q\right)_{n}}{\left(1-q\right)^{n}\prod_{j=1}^{s}\left(b_{j};q\right)_{n}}\\
&\times\, _{r}R_{s,q}\left(a_{1}q^{n},a_{2}q^{n},a_{i}q^{n};b_{1}q^{n},b_{2}q^{n},b_{j}q^{n};\alpha,n\alpha+\beta;q;q^{(1+s-r)n}z\right).
\end{split}
\end{equation}
\end{theorem}
\begin{proof}
By using the differential operator $\theta_{z,q}=zD_{z,q}$, we have
\begin{equation} \label{3.19}
\begin{split}
&\theta_{z,q}\;_{r}R_{s,q}=zD_{z,q}\;_{r}R_{s,q}=\sum_{\ell=0}^{\infty}\frac{1-q^{\ell}}{1-q}\frac{z^{\ell}\left(-1\right)^{\left(1+s-r\right)\ell}q^{\left(1+s-r\right)\binom{\ell}{2}}}{\left(q;q\right)_{\ell}\Gamma_{q}\left(\ell \alpha+\beta\right)}\frac{\prod_{i=1}^{r}\left(a_{i};q\right)_{\ell}}{\prod_{j=1}^{s}\left(b_{j};q\right)_{\ell}}\\
&=\frac{1}{1-q}\sum_{\ell=1}^{\infty}\frac{z^{\ell}(-1)^{\left(1+s-r\right)\ell}q^{\left(1+s-r\right)\binom{\ell}{2}}}{\left(q;q\right)_{\ell-1}\Gamma_{q}\left(\ell \alpha+\beta\right)}\frac{\prod_{i=1}^{r}(a_{i};q)_{\ell}}{\prod_{j=1}^{s}\left(b_{j};q\right)_{\ell}}\\
&=\frac{1}{1-q}\sum_{\ell=0}^{\infty}\frac{z^{\ell+1}\left(-1\right)^{\left(1+s-r\right)\left(\ell+1\right)}q^{\left(1+s-r\right)\binom{\ell+1}{2}}}{\left(q;q\right)_{\ell}\Gamma_{q}\left(\ell \alpha+\alpha+\beta\right)}\frac{\prod_{i=1}^{r}\left(a_{i};q\right)_{\ell+1}}{\prod_{j=1}^{s}\left(b_{j};q\right)_{\ell+1}}\\
&=\frac{z\left(-q\right)^{1+s-r}\prod_{i=1}^{r}\left(1-a_{i}\right)}{\left(1-q\right)\prod_{j=1}^{s}\left(1-b_{j}\right)}\, _{r}R_{s,q}\left(a_{1}q,a_{2}q,a_{i}q;b_{1}q,b_{2}q,b_{j}q;\alpha,\alpha+\beta;q;q^{1+s-r}z\right).
\end{split}
\end{equation}
Again, in general, we obtain	
\begin{equation*}
\begin{split}
\theta_{z,q}^{n}\;_{r}R_{s,q}=&\left(zD_{z,q}\right)^{n}\,_{r}R_{s,q}=\frac{z^{n}\left(-q\right)^{\left(1+s-r\right)n}\prod_{i=1}^{r}\left(a_{i};q\right)_{n}}{\left(1-q\right)^{n}\prod_{j=1}^{s}\left(b_{j};q\right)_{n}}\\
&\times\;_{r}R_{s,q}\left(a_{1}q^{n},a_{2}q^{n},a_{i}q^{n};b_{1}q^{n},b_{2}q^{n},b_{j}q^{n};\alpha,n\alpha+\beta;q;q^{\left(1+s-r\right)n}z\right).
\end{split}
\end{equation*}
\end{proof}
\begin{theorem}
The $\,_{r}R_{s,q}$-function has the following differential recurrence relations:
\begin{equation} \label{3.20}
	a_{i}\theta_{z,q}\;_{r}R_{s,q} = \frac{1-a_{i}}{1-q}\left(\,_{r}R_{s,q}(a_{i}q)-\;_{r}R_{s,q}\right)\;\;\; \left(i=1,2,\ldots,r\right),
	\end{equation}
and
\begin{equation} \label{3.21}
b_{j}q^{-1}\theta_{z,q}\,_{r}R_{s,q}=\frac{1-b_{j}q^{-1}}{1-q}\left(\,_{r}R_{s,q}(b_{j}q^{-1})-\;_{r}R_{s,q}\right)\;\;\; \left(j=1,2,\ldots,s\right).
	\end{equation}
\end{theorem}
\begin{proof}
By using \eqref{3.19}, the left-hand side of \eqref{3.20} can be written as
\begin{equation*}
	\begin{split}
	&\,_{r}R_{s,q}\left(a_{i}q\right)=\sum_{\ell=0}^{\infty}\frac{1-a_{i}q^{\ell}}{1-a_{i}}U_{\ell}=\sum_{\ell=0}^{\infty}\left[1+a_{i}\frac{1-q^{\ell}}{1-a_{i}}\right]U_{\ell}\\
	&=\;_{r}R_{s,q}+\frac{a_{i}\left(1-q\right)}{1-a_{i}}z\,D_{z,q}\,_{r}R_{s,q}=\, _{r}R_{s,q}+\frac{a_{i}\left(1-q\right)}{1-a_{i}}\theta_{z,q}\, _{r}R_{s,q}\;\;\; \left(i=1,2,\ldots,r\right),
	\end{split}
	\end{equation*}
and
\begin{equation*}
	\begin{split}
	&\, _{r}R_{s,q}\left(b_{j}q^{-1}\right)=\sum_{\ell=0}^{\infty}\frac{1-b_{j}q^{\ell-1}}{1-b_{j}q^{-1}}U_{\ell}=\sum_{\ell=0}^{\infty}\left[1+b_{j}q^{-1}\frac{1-q^{\ell}}{1-b_{j}q^{-1}}\right]U_{\ell}\\
	&=\sum_{\ell=0}^{\infty}U_{\ell}+\frac{b_{j}q^{-1}\left(1-q\right)}{1-b_{j}q^{-1}}\sum_{\ell=0}^{\infty}\frac{1-q^{\ell}}{1-q}U_{\ell}=\,_{r}R_{s,q}+\frac{b_{j}\left(1-q\right)}{q-b_{j}}\sum_{\ell=0}^{\infty}\left[\ell\right]_{q}U_{\ell}\\
	&=\,_{r}R_{s,q}+\frac{b_{j}(1-q)}{q-b_{j}}z\,D_{z,q}\,_{r}R_{s,q}=\,_{r}R_{s,q}+\frac{b_{j}\left(1-q\right)}{q-b_{j}}\theta_{z,q}\,_{r}R_{s,q}\;\;\;\left(j=1,2,\ldots,s\right).
	\end{split}
	\end{equation*}
Upon simplification, we obtain the results given in Eqs. \eqref{3.20} and \eqref{3.21}.	\end{proof}
\begin{theorem}
The $\,_{r}R_{s,q}$ function satisfies the following difference equations:
\begin{equation} \label{3.22}
	\left(D_{a_{i},q}\right)^{\kappa}\,_{r}R_{s,q} = \frac{(-1)^{\kappa}q^{\binom{\kappa}{2}}}{\left(a_{i};q\right)_{\kappa}}z^{\kappa}\left(D_{z,q}\right)^{\kappa}\,_{r}R_{s,q},
	\end{equation}
and
\begin{equation}   \label{3.23}
\left(D_{b_{j},q}\right)^{\kappa}\,_{r}R_{s,q}=\frac{z^{\kappa}}{\left(b_{j};q\right)_{\kappa}}\left(D_{z,q}\right)^{\kappa}z^{\kappa-1}\, _{r}R_{s,q}\left(b_{j}q^{\kappa}\right).
	\end{equation}
\end{theorem}
\begin{proof}
By using the $q$-difference operator $D_{a_{i},q}$, we have
\begin{equation}   \label{3.24}
	\begin{split}
	D_{a_{i},q}\;_{r}R_{s,q}=&\frac{\;_{r}R_{s,q}(a_{i}q)-\,_{r}R_{s,q}}{\left(q-1\right)a_{i}}\\
	=&\sum_{\ell=0}^{\infty}\frac{z^{\ell}(-1)^{\left(1+s-r\right)\ell}q^{\left(1+s-r\right)\binom{\ell}{2}}}{\left(q;q\right)_{\ell}\Gamma_{q}\left(\ell \alpha+\beta\right)}\,\frac{\prod_{i=1}^{r}\left(a_{i};q\right)_{\ell}}{\left(1-q\right)a_{i}\prod_{j=1}^{s}\left(b_{j};q\right)_{\ell}}\left(\frac{1-a_{i}q^{\ell}}{1-a_{i}}-1\right)\\
	=&-\frac{1}{1-a_{i}}\sum_{\ell=0}^{\infty}\frac{1-q^{\ell}}{1-q}\frac{z^{\ell}\left(-1\right)^{\left(1+s-r\right)\ell}q^{\left(1+s-r\right)\binom{\ell}{2}}}{\left(q;q\right)_{\ell}\Gamma_{q}\left(\ell \alpha+\beta\right)}\frac{\prod_{i=1}^{r}\left(a_{i};q\right)_{\ell}}{\prod_{j=1}^{s}\left(b_{j};q\right)_{\ell}}\\
	=&-\frac{z}{1-a_{i}}\,D_{z,q}\,_{r}R_{s,q}.
	\end{split}
	\end{equation}
By iterating \eqref{3.24} with respect to $\kappa$, we obtain the result in \eqref{3.22}. Similarly, the result in \eqref{3.23} can be proven.
\end{proof}
\begin{corollary}
The $\,_{r}R_{s,q}$ function satisfies the difference equations
\begin{equation} \label{3.25}
\left[\left(a_{i};q\right)_{\kappa}\left(D_{a_{i},q}\right)^{\kappa}-(a_{n};q)_{\kappa}\left(D_{a_{n},q}\right)^{\kappa}\right]\,_{r}R_{s,q}=0\;\;\; \left(i,n=1,2,...,r\right),
	\end{equation}
and
\begin{equation}   \label{3.26}
\left[\left(b_{j};q\right)_{\kappa}\left(D_{b_{j},q}\right)^{\kappa}-\left(b_{m};q\right)_{\kappa}\left(D_{b_{m},q}\right)^{\kappa}\right]\,_{r}R_{s,q}=0\;\;\; \left(j,m=1,2,...,s\right).
	\end{equation}
\end{corollary}
\begin{theorem}
The $\,_{r}R_{s,q}$ function has the following differential property:
\begin{equation}   \label{3.27}
	\begin{split}
	&\left(D_{z,q}\right)^{\kappa}\left[z^{\beta-1}\,_{r}R_{s,q}\left(a_{1},a_{2},\ldots,a_{r};b_{1},b_{2},\ldots,b_{s};\alpha,\beta;z^{\alpha}\right)\right]\\
	&=z^{\beta-\kappa-1}\,_{r}R_{s,q}\left(a_{1},a_{2},\ldots,a_{r};b_{1},b_{2},\ldots,b_{s}\alpha,\beta-\kappa;z^{\alpha}\right),
	\end{split}
	\end{equation}
	\begin{equation} \label{3.28}
	\begin{split}
	&D_{z,q}^{\kappa}\left[z^{a_{i}+\kappa-1}\,_{r}R_{s,q}\left(a_{1},a_{2},\ldots,a_{r};b_{1},b_{2},\ldots,b_{s};\alpha,\beta;z\right)\right]\\
	&=\frac{\Gamma_{q}\left(a_{i}+\kappa\right)}{\Gamma_{q}\left(a_{i}\right)}z^{a_{i}-1}\,_{r}R_{s,q}\left(a_{i}+\kappa,b_{j};\alpha,\beta;z\right),
	\end{split}
	\end{equation}
and
\begin{equation} \label{3.29}
	\begin{split}
	&D_{z,q}^{\kappa}\left[z^{b_{j}-1}\,_{r}R_{s,q}\left(a_{1},a_{2},\ldots,a_{r};b_{1},b_{2},\ldots,b_{s};\alpha,\beta;z\right)\right]\\
&=\frac{\Gamma_{q}\left(b_{j}\right)}{\Gamma_{q}\left(b_{j}-\kappa\right)}z^{b_{j}-1}\;_{r}R_{s,q}\left(a_{i},b_{j}-\kappa;\alpha,\beta;z\right).
	\end{split}
	\end{equation}
\end{theorem}
\begin{proof}
both sides of equation \eqref{2.1} by $z^{\beta-1}$ and differentiating term by term under the summation sign, we obtain
\begin{equation*}
	\begin{split}
	&\left(D_{z,q}\right)^{\kappa}\left[z^{\beta-1}\;_{r}R_{s,q}\left(a_{1},a_{2},\ldots,a_{r};b_{1},b_{2},\ldots,b_{s};\alpha,\beta;z^{\alpha}\right)\right]\\
	&=\sum_{\ell=0}^{\infty}\frac{(-1)^{\left(1+s-r\right)\left(\ell\right)}q^{(1+s-r)\binom{\ell}{2}}}{\left(q;q\right)_{\ell}\Gamma_{q}\left(\ell \alpha+\beta\right)}\frac{\prod_{i=1}^{r}\left(a_{i};q\right)_{\ell}}{\prod_{j=1}^{s}\left(b_{j};q\right)_{\ell}}\frac{\Gamma_{q}\left(\ell \alpha+\beta\right)}{\Gamma_{q}\left(\ell \alpha+\beta-\kappa\right)}z^{\ell\alpha+\beta-\kappa-1}\\
	&=z^{\beta-\kappa-1}\,_{r}R_{s,q}\left(a_{1},a_{2},\ldots,a_{r};b_{1},b_{2},\ldots,b_{s}\alpha,\beta-\kappa;z^{\alpha}\right).
	\end{split}
	\end{equation*}
	Similarly, the results in \eqref{3.28}–\eqref{3.29} can be proven.	
\end{proof}

\begin{theorem}
The $\,_{r}R_{s,q}$ function has the following differential properties:
\begin{equation} \label{3.30}
\begin{split}
D_{z,q}\left[z^{a_{1}}\,_{r}R_{s,q}\right]&=\frac{1-q^{a_{i}}}{1-q}z^{a_{i}-1}\;_{r}R_{s,q}+\frac{1-a_{i}}{\left(1-q\right)a_{i}}q^{a_{i}}z^{a_{i}-1}\left(\,_{r}R_{s,q}\left(a_{i}q\right)-\,_{r}R_{s,q}\right)\\
&=\frac{1-q^{a_{i}}}{1-q}z^{a_{i}-1}\,_{r}R_{s,q}+\frac{1-b_{j}q^{-1}}{\left(1-q\right)b_{j}q^{-1}}q^{a_{i}}z^{a_{i}-1}\left(\,_{r}R_{s,q}\left(b_{j}q^{-1}\right)-\,_{r}R_{s,q}\right),
	\end{split}
	\end{equation}
\begin{equation}   \label{3.31}
	\begin{split}
	D_{z,q}\left[z^{a_{1}}\;_{r}R_{s,q}\right]&=z^{a_{i}-1}\frac{1-a_{i}}{\left(1-q\right)a_{i}}\left(\,_{r}R_{s,q}\left(a_{i}q\right)-\,_{r}R_{s,q}\right)+\frac{1-q^{a_{i}}}{1-q}z^{a_{i}-1}\,_{r}R_{s,q}\left(qz\right),\\
	&=z^{a_{i}-1}\frac{1-b_{j}q^{-1}}{\left(1-q\right)b_{j}q^{-1}}\left(\,_{r}R_{s,q}\left(b_{j}q^{-1}\right)-\,_{r}R_{s,q}\right)+\frac{1-q^{a_{i}}}{1-q}z^{a_{i}-1}\,_{r}R_{s,q}\left(qz\right),
	\end{split}
	\end{equation}
\begin{equation} \label{3.32}
	\begin{split}
	D_{z,q}\left[z^{b_{j}+\ell-1}\,_{r}R_{s,q}\right]&=z^{b_{j}-2}\left[\left[b_{j}-1\right]_{q}\,_{r}R_{s,q}+
	\frac{\left(1-a_{i}\right)q^{b_{j}-1}}{\left(1-q\right)a_{i}}\left(\,_{r}R_{s,q}\left(a_{i}q\right)-\,_{r}R_{s,q}\right)\right]\\
&	=z^{b_{j}-2}\left[\left[b_{j}-1\right]_{q}\,_{r}R_{s,q}+
	\frac{\left(1-b_{j}q^{-1}\right)q^{b_{j}-1}}{\left(1-q\right)b_{j}q^{-1}}\left(\,_{r}R_{s,q}(b_{j}q^{-1})-\,_{r}R_{s,q}\right)\right]
	\end{split}
	\end{equation}
and
\begin{equation}  \label{3.33}
	\begin{split}
	D_{z,q}\left[z^{b_{j}+\ell-1}\,_{r}R_{s,q}\right]&=z^{b_{j}-2}\left[\frac{1-a_{i}}{\left(1-q\right)a_{i}}\left(\,_{r}R_{s,q}(a_{i}q)-\;_{r}R_{s,q}\right)+\left[b_{j}-1\right]_{q}\,_{r}R_{s,q}\left(qz\right)\right]\\
	&=z^{b_{j}-2}\left[\frac{1-b_{j}q^{-1}}{\left(1-q\right)b_{j}q^{-1}}\left(\,_{r}R_{s,q}\left(b_{j}q^{-1}\right)-\,_{r}R_{s,q}\right)+\left[b_{j}-1\right]_{q}\,_{r}R_{s,q}(qz)\right].
	\end{split}
	\end{equation}
\end{theorem}
\begin{proof}
By using the formula
\begin{equation} \label{3.34}
\left[\mu+\nu\right]_{q}=\left[\mu\right]_{q}+q^{\mu}\left[\nu\right]_{q}=\left[\nu\right]_{q}+q^{\nu}\left[\mu\right]_{q},
	\end{equation}
and
\begin{equation}  \label{3.35}
	\begin{split}
	D_{z,q}\left[z^{a_{i}+\ell}\right]=\left[a_{i}+\ell\right]_{q}z^{a_{i}+\ell-1}=\left[\left[a_{i}\right]_{q}+q^{a_{i}}\left[\ell\right]_{q}\right]z^{a_{i}+\ell-1}=\left[\left[\ell\right]_{q}+q^{\ell}\left[a_{i}\right]_{q}\right]z^{a_{i}+\ell-1}.
	\end{split}
	\end{equation}
Multiplying by $z^{a_{i}}$ in the equation \eqref{2.1} and differentiating term by term under the summation sign and using \eqref{3.34}-\eqref{3.35}, we obtain the result
\begin{equation*}
\begin{split}
	&D_{z,q}\left[z^{a_{i}}\,_{r}R_{s,q}\left(a_{1},a_{2},\ldots,a_{r};b_{1},b_{2},\ldots,b_{s};\alpha,\beta;z\right)\right]\\
	&=\sum_{\ell=0}^{\infty}\frac{\left[a_{i}+\ell\right]_{q}z^{a_{i}+\ell-1}(-1)^{\left(1+s-r\right)\left(\ell\right)}q^{\left(1+s-r\right)\binom{\ell}{2}}}{\left(q;q\right)_{\ell}\Gamma_{q}\left(\ell \alpha+\beta\right)}\frac{\prod_{i=1}^{r}\left(a_{i};q\right)_{\ell}}{\prod_{j=1}^{s}\left(b_{j};q\right)_{\ell}}\\
	&=\sum_{\ell=0}^{\infty}\left[\left[a_{i}\right]_{q}+q^{a_{i}}\left[\ell\right]_{q}\right]\frac{z^{a_{i}+\ell-1}(-1)^{\left(1+s-r\right)\left(\ell\right)}q^{\left(1+s-r\right)\binom{\ell}{2}}}{\left(q;q\right)_{\ell}\Gamma_{q}\left(\ell \alpha+\beta\right)}\frac{\prod_{i=1}^{r}\left(a_{i};q\right)_{\ell}}{\prod_{j=1}^{s}\left(b_{j};q\right)_{\ell}}\\
	&=\sum_{\ell=0}^{\infty}\frac{\left[a_{i}\right]_{q}z^{a_{i}+\ell-1}(-1)^{\left(1+s-r\right)\left(\ell\right)}q^{\left(1+s-r\right)\binom{\ell}{2}}}{\left(q;q\right)_{\ell}\Gamma_{q}\left(\ell \alpha+\beta\right)}\frac{\prod_{i=1}^{r}\left(a_{i};q\right)_{\ell}}{\prod_{j=1}^{s}\left(b_{j};q\right)_{\ell}}\\
	&+\sum_{\ell=0}^{\infty}\frac{q^{a_{i}}\left[\ell\right]_{q}z^{a_{i}+\ell-1}(-1)^{\left(1+s-r\right)\left(\ell\right)}q^{\left(1+s-r\right)\binom{\ell}{2}}}{\left(q;q\right)_{\ell}\Gamma_{q}\left(\ell \alpha+\beta\right)}\frac{\prod_{i=1}^{r}\left(a_{i};q\right)_{\ell}}{\prod_{j=1}^{s}\left(b_{j};q\right)_{\ell}}\\
	&=\left[a_{i}\right]_{q}z^{a_{i}-1}\,_{r}R_{s,q}\left(a_{1},a_{2},\ldots,a_{r};b_{1},b_{2},\ldots,b_{s}\alpha,\beta;z\right)\\
	&+q^{a_{i}}z^{a_{i}-1}\theta_{z,q}\;_{r}R_{s,q}\left(a_{1},a_{2},\ldots,a_{r};b_{1},b_{2},\ldots,b_{s}\alpha,\beta;z\right)\\
	&=\frac{1-q^{a_{i}}}{1-q}z^{a_{i}-1}\,_{r}R_{s,q}+\frac{1-a_{i}}{\left(1-q\right)a_{i}}q^{a_{i}}z^{a_{i}-1}\left(\,_{r}R_{s,q}\left(a_{i}q\right)-\,_{r}R_{s,q}\right).
	\end{split}
	\end{equation*}
The proofs of \eqref{3.31}–\eqref{3.33} follow similarly to the proof above.
\end{proof}

\section{Integral representations for the $\,_{r}R_{s,q}$ function} \label{sec 4}
Here, we establish the integral representations of the $\;_{r}R_{s,q}$ function in terms of other well-known special functions of fractional calculus, and these can be extended to derive the following results:
\begin{theorem}
The $\,_{r}R_{s,q}$ function have the following integral representation:
\begin{equation} \label{4.1}
	\begin{split}
	&\;_{r}R_{s,q}\left(a_{1},a_{2},\ldots ,a_{r};b_{1},b_{2},\ldots ,b_{s};\alpha,\beta,z\right)\\
	&=\frac{\Gamma_{q} \left(b_{j}\right)}{\Gamma_{q} \left(a_{i}\right)\Gamma_{q} \left(b_{j}-a_{i}\right)} \int_{0}^{1}t^{a_{i}-1}\left(1-qt\right)^{b_{j}-a_{i}-1} \\
	&\times\;_{r-1}R_{s-1,q}\left(\begin{array}{cc}
	a_{1},\ldots,a_{i-1},a_{i+1}\ldots,a_{r} \\
	b_{1},\ldots,b_{j-1},b_{j+1}\ldots,b_{s}
	\end{array} ; \alpha,\beta,zt\right).
	\end{split}
	\end{equation}
\end{theorem}
\begin{proof}
By using integral definition of $q$-beta function \eqref{1.10}, we get
\begin{equation}   \label{4.3}
\frac{\left(a_{i};q\right)_{\ell}}{\left(b_{j};q\right)_{\ell}} = \frac{\Gamma_{q} \left( b_{j}\right)}{\Gamma_{q} \left(a_{i}\right)\Gamma_{q} \left(b_{j}-a_{i}\right)} \int_{0}^{1}t^{a_{i}+\ell-1}\left(1-qt\right)^{b_{j}-a_{i}-1}\mathrm{d}_{q}t,
\end{equation}
by applying equations \eqref{4.3} and \eqref{2.1}, we obtain \eqref{4.1}.
\end{proof}
\begin{corollary}
The following $q$-integral for $_{r}R_{s,q}$ holds true:
\begin{equation} \label{4.4}
	\begin{split}
	&\,_{r}R_{s,q}\left(a_{1},a_{2},\ldots,a_{r};b_{1},b_{2},\ldots,b_{s};\alpha,\beta;z\right)\\
	&=\prod_{i=1}^{r}\prod_{j=1}^{s}\frac{1}{B_{q}\left(a_{i},b_{j}-a_{i}\right)}\int_{0}^{1}t^{a_{i}-1}\left(qt;q\right)_{b_{j}-a_{i}-1}\,_{r}R_{s,q}\left(0,0,\ldots,0;0,0,\ldots,0;\alpha,\beta;zt\right)\mathrm{d}_{q}t.
	\end{split}
	\end{equation}
\end{corollary}
\begin{proof} By using formula \eqref{4.3}, we obtain
\begin{equation*}
	\begin{split}
	&\,_{r}R_{s,q}\left(a_{1},a_{2},\ldots,a_{r};b_{1},b_{2},\ldots,b_{s};\alpha,\beta;z\right)\\
	&=\sum_{\ell=0}^{\infty}\frac{z^{\ell}(-1)^{\left(1+s-r\right)\left(\ell\right)}q^{\left(1+s-r\right)\binom{\ell}{2}}}{\left(q;q\right)_{\ell}\Gamma_{q}\left(\ell \alpha+\beta\right)}\prod_{i=1}^{r}\frac{\Gamma_{q}\left(a_{i}+\ell\right)\left(1-q\right)^{\ell}}{\Gamma_{q}\left(a_{i}\right)}\prod_{j=1}^{s}\frac{\Gamma_{q}\left(b_{j}\right)}{\Gamma_{q}\left(b_{j}+\ell\right)\left(1-q\right)^{\ell}}\\
	&=\sum_{\ell=0}^{\infty}\frac{z^{\ell}(-1)^{\left(1+s-r\right)\left(\ell\right)}q^{\left(1+s-r\right)\binom{\ell}{2}}}{\left(q;q\right)_{\ell}\Gamma_{q}\left(\ell \alpha+\beta\right)}\prod_{i=1}^{r}\prod_{j=1}^{s}\frac{B_{q}\left(a_{i}+\ell,b_{j}-a_{i}\right)}{B_{q}\left(a_{i},b_{j}-a_{i}\right)}\\
	&=\prod_{i=1}^{r}\prod_{j=1}^{s}\frac{1}{B_{q}\left(a_{i},b_{j}-a_{i}\right)}\sum_{\ell=0}^{\infty}\frac{z^{\ell}\left(-1\right)^{\left(1+s-r\right)\left(\ell\right)}q^{\left(1+s-r\right)\binom{\ell}{2}}}{\left(q;q\right)_{\ell}\Gamma_{q}\left(\ell \alpha+\beta\right)}\int_{0}^{1}\left(zt\right)^{\ell}t^{a_{i}-1}\left(qt;q\right)_{b_{j}-a_{i}-1}\mathrm{d}_{q}t\\
	&=\prod_{i=1}^{r}\prod_{j=1}^{s}\frac{1}{B_{q}\left(a_{i},b_{j}-a_{i}\right)}\int_{0}^{1}t^{a_{i}-1}\left(qt;q\right)_{b_{j}-a_{i}-1}\,_{r}R_{s,q}\left(0,0,\ldots,0;0,0,\ldots,0;\alpha,\beta;zt\right)\mathrm{d}_{q}t.
	\end{split}
	\end{equation*}
\end{proof}

\section{$q$-Laplace, $q$-Sumudu, and $q$-Natural Transforms}  \label{sec 5}
In this section, we evaluate the $q$-Laplace, $q$-Sumudu and $q$-Natural transforms for the $q$-analogue of $\,_{r}R_{s,q}$ function.
\begin{theorem}
The $q$-Laplace transform of the $_{r}R_{s,q}$ function is determined by
\begin{equation}   \label{5.1}
	\begin{split}
	&\mathfrak{L}_{q}\left[t^{\beta-1}\,_{r}R_{s,q}\left(a_{1},a_{2},\ldots,a_{r};b_{1},b_{2},\ldots,b_{s};\alpha,\beta,zt^{\alpha}\right);s\right]\\
	&=s^{-\beta}\left(1-q\right)^{\beta-1}\,_{r}\Phi_{s}\left(a_{1},a_{2},\ldots,a_{r};b_{1},b_{2},\ldots,b_{s};z\left(1-q\right)^{\alpha}s^{-\alpha}\right),
	\end{split}
	\end{equation}
where $\mathfrak{L}_{q}[f(t);s]$ is the $q$-Laplace transform, defined as
\begin{equation*}
	\mathfrak{L}_{q}\left[f(t);s\right]=\frac{1}{1-q}\int_{0}^{\frac{1}{s}}E_{q}^{qst}f(t)\,\mathrm{d}_{q}t=F(s)\;\;\; \left(s\in\mathbb{C}\right).
	\end{equation*}
\end{theorem}
\begin{proof}
Using the $q$-Euler's integral, we get
\begin{equation}   \label{5.2}
\mathfrak{L_{q}}\left[t^{\ell \alpha+\beta-1};s\right]=\frac{1}{1-q}\int_{0}^{\frac{1}{s}}E_{q}^{qst}\,t^{\ell \alpha+\beta-1}\mathrm{d}_{q}t=\frac{\Gamma_{q}(\ell \alpha+\beta)\left(1-q\right)^{\ell \alpha+\beta-1}}{s^{\ell \alpha+\beta}},
	\end{equation}
by using the above equation \eqref{5.2}, this yields the right-hand side of \eqref{5.1}.
\end{proof}
\begin{theorem}
The $q$-Sumudu transform of $_{r}R_{s,q}$ function is determined by
\begin{equation}   \label{5.3}
	\begin{split}
	&\mathfrak{S}_{q}\left[t^{\beta-1}\,_{r}R_{s,q}\left(a_{1},a_{2},\ldots,a_{r};b_{1},b_{2},\ldots,b_{s};\alpha,\beta,zt^{\alpha}\right);s\right]\\
	&=s^{\beta-1}\left(1-q\right)^{\beta-1}\,_{r}\Phi_{s}\left(a_{1},a_{2},\ldots,a_{r};b_{1},b_{2},\ldots,b_{s};z(1-q)^{\alpha}s^{\alpha}\right),
	\end{split}
	\end{equation}
where $\mathfrak{S}_{q}[f(t);s]$ is $q$-Sumudu transform, defined as
\begin{equation*}
	\begin{split}
	\mathfrak{S}_{q}\left\{f(t);s\right\}&=\frac{1}{\left(1-q\right)s}\int_{0}^{s}f(t)\,E_{q}\left(\frac{qt}{s}\right)\mathrm{d}_{q}t\\
	&=\frac{1}{1-q}\int_{0}^{1}f(st)E_{q}(qt)\mathrm{d}_{q}t\\
	\mathfrak{S}_{q}\left\{t^{\alpha-1};s\right\}&=s^{\alpha-1}\left(1-q\right)^{\alpha-1}\Gamma_{q}\left(\alpha\right)\;\; \left(s\in \left(-\tau_{1},\tau_{2}\right)\right),
	\end{split}
	\end{equation*}
over the set of functions
\begin{equation*}
A=\left\{f(t)\exists M,\tau_{1}, \tau_{2}>0,\,\left|f(t)\right|<ME_{q}\left(\frac{\left|t\right|}{\tau_{j}}\right),\,t\in (-1)^{j}\times  \left[0,\infty\right[\right\}.
	\end{equation*}
\end{theorem}
\begin{proof}
By applying the $q$-Sumudu transform to \eqref{2.1}, we obtain \eqref{5.3}.	
\end{proof}
\begin{theorem}
The $q$-Sumudu transform of the $\,_{r}R_{s,q}$ function is given as
\begin{equation}     \label{5.4}
\begin{split}
	&\mathfrak{S}_{q}\left\{\;_{r}R_{s,q}\left(a_{1},a_{2},\ldots,a_{r};b_{1},b_{2},\ldots,b_{s};\alpha,\beta;\alpha,\beta;zt\right);s\right\}\\
&=\, _{r}M_{s,q}\left(a_{1},a_{2},\ldots,a_{r};b_{1},b_{2},\ldots,b_{s};\alpha,\beta;\alpha,\beta;z s\right).
	\end{split}
	\end{equation}
\end{theorem}
\begin{proof}
Applying the $q$-analogues of the $q$-Sumudu transform in \eqref{2.1}, we get
\begin{equation*}
\begin{split}
&\mathfrak{S}_{q}\left\{\,_{r}R_{s,q}\left(a_{1},a_{2},\ldots,a_{r};b_{1},b_{2},\ldots,b_{s};\alpha,\beta;\alpha,\beta;zt\right);s\right\}\\
&=\sum_{\ell=0}^{\infty}\frac{z^{\ell}(-1)^{\left(1+s-r\right)\left(\ell\right)}q^{\left(1+s-r\right)\binom{\ell}{2}}}{\left(q;q\right)_{\ell}\Gamma_{q}\left(\ell \alpha+\beta\right)}\frac{\prod_{i=1}^{r}\left(a_{i};q\right)_{\ell}}{\prod_{j=1}^{s}\left(b_{j};q\right)_{\ell}}\mathfrak{S}_{q}\left\{t^{\ell};s\right\}\\
& =\sum_{\ell=0}^{\infty}\frac{z^{\ell}\left(-1\right)^{\left(1+s-r\right)\left(\ell\right)}q^{\left(1+s-r\right)\binom{\ell}{2}}}{\left(q;q\right)_{\ell}\Gamma_{q}\left(\ell \alpha+\beta\right)}\frac{\prod_{i=1}^{r}(a_{i};q)_{\ell}}{\prod_{j=1}^{s}\left(b_{j};q\right)_{\ell}}s^{\ell}\left(1-q\right)^{\ell}\Gamma_{q}\left(\ell+1\right)\\
&=\sum_{\ell=0}^{\infty}\frac{z^{\ell}\left(-1\right)^{\left(1+s-r\right)\left(\ell\right)}q^{\left(1+s-r\right)\binom{\ell}{2}}}{\Gamma_{q}\left(\ell \alpha+\beta\right)}\frac{\prod_{i=1}^{r}\left(a_{i};q\right)_{\ell}}{\prod_{j=1}^{s}\left(b_{j};q\right)_{\ell}}s^{\ell}\\
&=\,_{r}M_{s,q}\left(a_{1},a_{2},\ldots,a_{r};b_{1},b_{2},\ldots,b_{s};\alpha,\beta;\alpha,\beta;z s\right).
\end{split}
\end{equation*}
\end{proof}

\begin{theorem}
Let $\mu$ be a positive real number. Then, the $q$-Natural transform of the $\;_{r}R_{s,q}$ function is given by
\begin{equation}   \label{5.5}
	\begin{split}
	&\mathbf{N}_{q}\left\{t^{\mu-1}\,_{r}R_{s,q}\left(\begin{array}{cc}
	a_{1},\ldots,a_{i-1},a_{i+1}\ldots,a_{r} \\
	b_{1},\ldots,b_{j-1},b_{j+1}\ldots ,b_{s}
	\end{array};\, \alpha,\beta,zt\right);\,(u;s)\right\}\\
	&=\frac{\Gamma_{q}(\mu)\left(1-q\right)^{\mu-1}u^{\mu-1}}{s^{\mu}}
	\,_{r+1}R_{s+1,q}\left(\begin{array}{cc}
a_{1},\ldots,a_{i-1},a_{i+1}\ldots ,a_{r}, q^{\mu} \\
b_{1},\ldots ,b_{j-1},b_{j+1}\ldots ,b_{s}, 0 \end{array}; \, \alpha,\beta,\frac{zu}{s}\right).
	\end{split}
	\end{equation}
\end{theorem}
\begin{proof}
By applying the identity of the gamma function, we get	
\begin{equation}   \label{5.6}
\begin{split}
\Gamma_{q}\left(\mu+\ell\right)=\frac{\left(q^{\mu};q\right)_{\ell}}{\left(1-q\right)^{\ell}}\Gamma_{q}\left(\mu\right) = \frac{\left(q^{\mu};q\right)_{\ell}}{\left(1-q\right)^{\ell}}\Gamma_{q}\left(\mu\right),
	\end{split}
	\end{equation}
By applying the  $q$-Natural transform to (\ref{2.1}) and (\ref{5.6}), we obtain
\begin{equation*}
	\begin{split}
&\mathbf{N}_{q}\left\{t^{\mu-1}\,_{r}R_{s,q}\left(\begin{array}{cc}
	a_{1},\ldots,a_{i-1},a_{i+1}\ldots ,a_{r} \\
	b_{1},\ldots,b_{j-1},b_{j+1}\ldots,b_{s}
	\end{array};\, \alpha,\beta,t\right);\, \left(u;s\right)\right\}\\
&=\sum_{\ell=0}^{\infty}\frac{z^{\ell}(-1)^{\left(1+s-r\right)\left(\ell\right)}q^{\left(1+s-r\right)\binom{\ell}{2}}}{\left(q;q\right)_{\ell}\Gamma_{q}\left(\ell \alpha+\beta\right)}\frac{\prod_{i=1}^{r}\left(a_{i};q\right)_{\ell}}{\prod_{j=1}^{s}\left(b_{j};q\right)_{\ell}}\mathbf{N}_{q}\left\{t^{\mu+\ell-1};\left(u;s\right)\right\}\\
&=\sum_{\ell=0}^{\infty}\frac{z^{\ell}(-1)^{\left(1+s-r\right)\left(\ell\right)}q^{\left(1+s-r\right)\binom{\ell}{2}}}{\left(q;q\right)_{\ell}\Gamma_{q}\left(\ell \alpha+\beta\right)}\frac{\prod_{i=1}^{r}\left(a_{i};q\right)_{\ell}}{\prod_{j=1}^{s}\left(b_{j};q\right)_{\ell}}\frac{u^{\mu+\ell-1}}{s^{\mu+\ell}}\left(1-q\right)^{\mu+\ell-1}\Gamma_{q}\left(\mu+\ell\right)\\
&=\sum_{\ell=0}^{\infty}\frac{z^{\ell}(-1)^{(1+s-r)\left(\ell\right)}q^{\left(1+s-r\right)\binom{\ell}{2}}}{\Gamma_{q}\left(\ell \alpha+\beta\right)}\frac{\prod_{i=1}^{r}\left(a_{i};q\right)_{\ell}}{\prod_{j=1}^{s}\left(b_{j};q\right)_{\ell}}\frac{u^{\mu+\ell-1}}{s^{\mu+\ell}}\frac{\left(q^{\mu};q\right)_{\ell}\left(1-q\right)^{\mu+\ell-1}}{\left(1-q\right)^{\ell}}\Gamma_{q}\left(\mu\right)\\
&=\frac{\Gamma_{q}\left(\mu\right)\left(1-q\right)^{\mu-1}u^{\mu-1}}{s^{\mu}}
	\,_{r+1}R_{s+1,q}\left(\begin{array}{cc}
	a_{1},\ldots,a_{i-1},a_{i+1}\ldots,a_{r}, q^{\mu} \\
	b_{1},\ldots,b_{j-1},b_{j+1}\ldots,b_{s}, 0
\end{array};\, \alpha,\beta,\frac{zu}{s}\right).
\end{split}
\end{equation*}
\end{proof}

\section{Fractional $q$-Integrals and $q$-Derivatives of $\;_{r}R_{s,q}$ series} \label{sec 6}
In this section, we derive and evaluate the fractional $q$-derivatives and $q$-integrals operators and properties of $\,_{r}R_{s,q}$ function by using fractional $q$-calculus.
\begin{theorem}
Let $\mu$ and be $\mathbf{I}_{b^{-},q}^{\alpha}$ the left-sided and right-sided operators of Riemann-Liouville fractional $q$-integral \eqref{1.6}-\eqref{1.7}. Then, the following formulas hold:
\begin{equation} \label{6.1}
\begin{split}
&\mathbf{I}_{0^{+},q}^{\mu}\left(t^{\beta-1}\,_{r}R_{s,q}\left(a_{1},a_{2},\ldots,a_{r};b_{1},b_{2},\ldots,b_{s};\alpha,\beta;\alpha,\beta;zt^{\alpha}\right)\right)\left(x\right)\\
&=x^{\mu+\beta-1}\,_{r}R_{s,q}\left(a_{1},a_{2},\ldots,a_{r};b_{1},b_{2},\ldots,b_{s};\alpha,\beta;\alpha,\beta+\mu;zx^{\alpha}\right),
	\end{split}
	\end{equation}
and
\begin{equation} \label{6.2}
	\begin{split}
&	\mathbf{I}_{0^{-},q}^{\mu}\left(t^{-\mu-\beta}\,_{r}R_{s,q}\left(a_{1},a_{2},\ldots,a_{r};b_{1},b_{2},\ldots,b_{s};\alpha,\beta;\alpha,\beta;zt^{-\alpha}\right)\right)\left(x\right)\\
&=\frac{x^{-\beta}}{q}\;_{r}R_{s,q}\left(a_{1},a_{2},\ldots,a_{r};b_{1},b_{2},\ldots,b_{s};\alpha,\beta;\alpha,\beta+\mu;zx^{-\alpha}\right).
	\end{split}
	\end{equation}
\end{theorem}
\begin{proof}
By applying the definition of the $\,_{r}R_{s,q}$-series in \eqref{2.1} and the fractional $q$-integral formula in \eqref{1.7}, we obtain
\begin{equation} \label{6.3}
	\begin{split}
	&\mathbf{I}_{0^{+},q}^{\mu}\left(t^{\beta-1}\,_{r}R_{s,q}\left(a_{1},a_{2},\ldots,a_{r};b_{1},b_{2},\ldots,b_{s};\alpha,\beta;\alpha,\beta;xt^{\alpha}\right)\right)(z)\\
&=\sum_{\ell=0}^{\infty}\frac{z^{\ell}\left(-1\right)^{\left(1+s-r\right)\left(\ell\right)}q^{\left(1+s-r\right)\binom{\ell}{2}}}{\left(q;q\right)_{\ell}\Gamma_{q}\left(\ell \alpha+\beta\right)}\frac{\prod_{i=1}^{r}\left(a_{i};q\right)_{\ell}}{\prod_{j=1}^{s}\left(b_{j};q\right)_{\ell}}\frac{x^{\mu-1}}{\Gamma_{q}\left(\mu\right)}\int_{0}^{x}\left(\frac{tq}{x};q\right)_{\mu-1}t^{\ell \alpha+\beta-1}\mathrm{d}_{q}t\\
&	=\sum_{\ell=0}^{\infty}\frac{z^{\ell}\left(-1\right)^{\left(1+s-r\right)\left(\ell\right)}q^{\left(1+s-r\right)\binom{\ell}{2}}}{\left(q;q\right)_{\ell}\Gamma_{q}\left(\ell \alpha+\beta\right)}\frac{\prod_{i=1}^{r}\left(a_{i};q\right)_{\ell}}{\prod_{j=1}^{s}\left(b_{j};q\right)_{\ell}}\frac{x^{\mu-1}}{\Gamma_{q}\left(\mu\right)}\int_{0}^{x}\frac{\left(\frac{tq}{x};q\right)_{\infty}}{\left(\frac{tq^{\mu}}{x};q\right)_{\infty}}t^{\ell \alpha+\beta-1}\mathrm{d}_{q}t,
	\end{split}
	\end{equation}
by substituting $t=xu$, $\mathrm{d}_{q}t=x\mathrm{d}_{q}u$, we get
\begin{equation} \label{6.4}
	\begin{split}
	\frac{x^{\mu-1}}{\Gamma_{q}(\mu)}\int_{0}^{x}\frac{\left(\frac{tq}{z};q\right)_{\infty}}{\left(\frac{tq^{\mu}}{z};q\right)_{\infty}}\,t^{\ell \alpha+\beta-1}\mathrm{d}_{q}t = \frac{x^{\mu+\ell \alpha+\beta-1}\Gamma_{q}\left(\ell \alpha+\beta\right)}{\Gamma_{q}\left(\ell \alpha+\beta+\mu\right)}.
	\end{split}
\end{equation}
By substituting \eqref{6.4} into \eqref{6.3} and simplifying, we obtain \eqref{6.1}. The proof of \eqref{6.2} follows similarly to the proof above.	
\end{proof}
\begin{remark}
Readers can refer to recent works \cite{GuptaShaktawatKumar2016,GuptaShaktawatKumar2017,KhanRamaniSutharKumar2018,Kumar2016,KumarKumar2014} on fractional calculus operators.
\end{remark}
\begin{theorem}
Let $\mu$ be given, along with the right-sided operator of the Riemann-Liouville $q$-fractional integral in \eqref{1.18}. Then, the following formulas hold:
\begin{equation}  \label{6.10}
\begin{split}
&\mathbf{D}_{0^{+},q}^{\mu}\left(t^{\beta-1}\,_{r}R_{s,q}\left(a_{1},a_{2},\ldots,a_{r};b_{1},b_{2},\ldots,b_{s};\alpha,\beta;\alpha,\beta;zt^{\alpha}\right)\right)(x)\\
&=x^{\beta-\mu-1}\,_{r}R_{s,q}\left(a_{1},a_{2},\ldots,a_{r};b_{1},b_{2},\ldots,b_{s};\alpha,\beta;\alpha,\beta-\mu;zx^{\alpha}\right).
	\end{split}
\end{equation}
\end{theorem}
\begin{proof}
By using \eqref{2.1} and $q$-fractional $q$-integral formula \eqref{1.18}, we obtain
\begin{equation} \label{6.12}
\begin{split}
	&\mathbf{D}_{0^{+},q}^{\mu}\left(t^{\beta-1}\,_{r}R_{s,q}\left(a_{1},a_{2},\ldots,a_{r};b_{1},b_{2},\ldots,b_{s};\alpha,\beta;\alpha,\beta;zt^{\alpha}\right)\right)(x)\\
&=\sum_{\ell=0}^{\infty}\frac{z^{\ell}(-1)^{\left(1+s-r\right)\left(\ell\right)}q^{\left(1+s-r\right)\binom{\ell}{2}}}{\left(q;q\right)_{\ell}\Gamma_{q}\left(\ell \alpha+\beta\right)}\frac{\prod_{i=1}^{r}\left(a_{i};q\right)_{\ell}}{\prod_{j=1}^{s}\left(b_{j};q\right)_{\ell}}\mathbf{D}_{0^{+},q}^{\mu}t^{\ell \alpha+\beta-1} = D^{n}_{q}\,\mathbf{I}_{0^{+},q}^{n-\mu}t^{\ell \alpha+\beta-1}\\
&=\sum_{\ell=0}^{\infty}\frac{z^{\ell}(-1)^{(1+s-r)(\ell)}q^{(1+s-r)\binom{\ell}{2}}}{\left(q;q\right)_{\ell}\Gamma_{q}\left(\ell \alpha+\beta\right)}\frac{\prod_{i=1}^{r}\left(a_{i};q\right)_{\ell}}{\prod_{j=1}^{s}\left(b_{j};q\right)_{\ell}}\frac{1}{\Gamma_{q}\left(n-\mu\right)}\\
&\times\;D^{n}_{q}\int_{0}^{x}x^{n-\mu-1}\left(\frac{qt}{x};q\right)_{n-\mu-1}t^{\ell \alpha+\beta-1}\mathrm{d}_{q}t\\
&=\sum_{\ell=0}^{\infty}\frac{z^{\ell}(-1)^{(1+s-r)(\ell)}\,q^{(1+s-r)\binom{\ell}{2}}}{\left(q;q\right)_{\ell}\Gamma_{q}\left(\ell \alpha+\beta\right)}\,\frac{\prod_{i=1}^{r}\left(a_{i};q\right)_{\ell}}{\prod_{j=1}^{s}\left(b_{j};q\right)_{\ell}}\frac{1}{\Gamma_{q}\left(n-\mu\right)}\\
& \times\; D^{n}_{q}x^{n-\mu-1}\int_{0}^{x}\frac{\left(\frac{qt}{x};q\right)_{\infty}}{\left(\frac{q^{n-\mu}t}{x};q\right)_{\infty}}\,t^{\ell \alpha+\beta-1}\,\mathrm{d}_{q}t,
	\end{split}
	\end{equation}
by substituting $t=xu$, $d_{q}t=xd_{q}u$, we get
\begin{equation}   \label{6.13}
\frac{1}{\Gamma_{q}\left(n-\mu\right)}\,D^{n}_{q}x^{n-\mu+\ell \alpha+\beta-1}\int_{0}^{1}\frac{\left(qu;q\right)_{\infty}}{\left(q^{n-\mu}u;q\right)_{\infty}}u^{\ell \alpha+\beta-1}\,\mathrm{d}_{q}u = \frac{x^{\beta-\mu-1}\Gamma_{q}\left(\ell \alpha+\beta\right)}{\Gamma_{q}\left(\ell \alpha+\beta-\mu\right)}x^{\ell \alpha}.	\end{equation}
By applying equations \eqref{6.13} and \eqref{6.12} and simplifying, the right-hand side (R.H.S.) of the above equation reduces to \eqref{6.10}.	
\end{proof}
\begin{theorem}
The following Kober fractional $q$-integral and $q$-derivative operators for the $_{r}R_{s,q}$ function hold:
\begin{equation}\label{6.19}
\begin{split}
&I_{q}^{\nu,\mu}\left[\,_{r}R_{s,q}\left(a_{1},a_{2},\ldots,a_{r};b_{1},b_{2},\ldots,b_{s};\alpha,\beta;z\right)\right]\\
&=\frac{\Gamma_{q}\left(\nu+1\right)}{\Gamma_{q}\left(\nu+\mu+1\right)}\,_{r+1}R_{s+1,q}\left(a_{1},a_{2},\ldots,a_{r},q^{\nu+1};b_{1},b_{2},\ldots,b_{s},q^{\nu+\mu+1};\alpha,\beta;z\right),
\end{split}
\end{equation}
and
\begin{equation} \label{6.20}
	\begin{split}
	&D_{q}^{\nu,\mu}\left[\,_{r}R_{s,q}\left(a_{1},a_{2},\ldots,a_{r};b_{1},b_{2},\ldots,b_{s};\alpha,\beta;z\right)\right]\\
	&=\frac{\Gamma_{q}\left(\nu+\mu+1\right)}{\Gamma_{q}\left(\nu+1\right)}\,_{r+1}R_{s+1,q}\left(a_{1},a_{2},\ldots,a_{r},q^{\nu+\mu+1};b_{1},b_{2},\ldots,b_{s},q^{\nu+1};\alpha,\beta;z\right).
	\end{split}
	\end{equation}
\end{theorem}
\begin{proof}
By applying the definition in \eqref{2.1} and the Kober fractional $q$-integral formula in \eqref{1.21}, we have
	\begin{equation*}
	\begin{split}
&I_{q}^{\nu,\mu}\left[\,_{r}R_{s,q}\left(a_{1},a_{2},\ldots,a_{r};b_{1},b_{2},\ldots,b_{s};\alpha,\beta;z\right)\right]\\
&=\sum_{\ell=0}^{\infty}\frac{c^{\ell}(-1)^{(1+s-r)(\ell)}q^{(1+s-r)\binom{\ell}{2}}}{\left(q;q\right)_{\ell}\Gamma_{q}\left(\ell \alpha+\beta\right)}\frac{\prod_{i=1}^{r}\left(a_{i};q\right)_{\ell}}{\prod_{j=1}^{s}\left(b_{j};q\right)_{\ell}}\,I_{q}^{\nu,\mu}z^{\ell}\\
&=\sum_{\ell=0}^{\infty}\frac{(-1)^{(1+s-r)(\ell)}q^{(1+s-r)\binom{\ell}{2}}}{\left(q;q\right)_{\ell}\Gamma_{q}\left(\ell \alpha+\beta\right)}\frac{\prod_{i=1}^{r}\left(a_{i};q\right)_{\ell}}{\prod_{j=1}^{s}\left(b_{j};q\right)_{\ell}}\frac{\Gamma_{q}\left(\nu+\ell+1\right)}{\Gamma_{q}\left(\nu+\mu+\ell+1\right)}z^{\ell}\\
&=\frac{\Gamma_{q}\left(\nu+1\right)}{\Gamma_{q}\left(\nu+\mu+1\right)}\sum_{\ell=0}^{\infty}\frac{(-1)^{(1+s-r)(\ell)}q^{(1+s-r)\binom{\ell}{2}}}{\left(q;q\right)_{\ell}\Gamma_{q}\left(\ell \alpha+\beta\right)}\frac{\prod_{i=1}^{r}\left(a_{i};q\right)_{\ell}}{\prod_{j=1}^{s}\left(b_{j};q\right)_{\ell}}\frac{\left(q^{\nu+1};q\right)_{\ell}}{\left(q^{\nu+\mu+1};q\right)_{\ell}}z^{\ell}\\
&=\frac{\Gamma_{q}\left(\nu+1\right)}{\Gamma_{q}\left(\nu+\mu+1\right)}\,_{r+1}R_{s+1,q}\left(a_{1},a_{2},\ldots,a_{r},q^{\nu+1};b_{1},b_{2},\ldots,b_{s},q^{\nu+\mu+1};\alpha,\beta;z\right).
\end{split}
\end{equation*}
Similarly, we obtain the result (\ref{6.20}).
\end{proof}
\begin{theorem}
The following assumption holds true:
\begin{equation} \label{6.21}
	\begin{split}
&K_{-,q}^{\nu,\mu}\left[\,_{r}R_{s,q}\left(a_{1},a_{2},\ldots,a_{r};b_{1},b_{2},\ldots,b_{s};\alpha,\beta;\frac{1}{z}\right)\right]\\
&	=\frac{\Gamma_{q}\left(\nu\right)}{\Gamma_{q}\left(\nu+\mu\right)}\,_{r+1}R_{s+1,q}\left(a_{1},a_{2},\ldots,a_{r},q^{\nu};b_{1},b_{2},\ldots,b_{s},q^{\nu+\mu};\alpha,\beta;\frac{1}{z}q^{\mu}\right),
\end{split}
\end{equation}
and
\begin{equation} \label{6.22}
\begin{split}
&D_{q}^{\nu,\mu}\left[\,_{r}R_{s,q}\left(a_{1},a_{2},\ldots,a_{r};b_{1},b_{2},\ldots,b_{s};\alpha,\beta;\frac{1}{z}\right)\right]\\
&=\frac{\Gamma_{q}\left(\nu+\mu\right)}{\Gamma_{q}(\nu)}\,_{r+1}R_{s+1,q}\left(a_{1},a_{2},\ldots,a_{r},q^{\nu+\mu};b_{1},b_{2},\ldots,b_{s},q^{\nu};\alpha,\beta;\frac{1}{z}q^{-\mu}\right).
\end{split}
\end{equation}
\end{theorem}
\begin{proof}
To prove \eqref{6.21}, we proceed by applying \eqref{1.22}, then we arrive at		\begin{equation*}
	\begin{split}
	&K_{-,q}^{\nu,\mu}\left[\,_{r}R_{s,q}\left(a_{1},a_{2},\ldots,a_{r};b_{1},b_{2},\ldots,b_{s};\alpha,\beta;\frac{1}{z}\right)\right]\\
&=\sum_{\ell=0}^{\infty}\frac{c^{\ell}(-1)^{(1+s-r)(\ell)}q^{(1+s-r)\binom{\ell}{2}}}{\left(q;q\right)_{\ell}\Gamma_{q}\left(\ell \alpha+\beta\right)}\frac{\prod_{i=1}^{r}\left(a_{i};q\right)_{\ell}}{\prod_{j=1}^{s}\left(b_{j};q\right)_{\ell}}\,I_{q}^{\nu,\mu}z^{-\ell}\\
&=\sum_{\ell=0}^{\infty}\frac{(-1)^{(1+s-r)(\ell)}q^{(1+s-r)\binom{\ell}{2}}}{\left(q;q\right)_{\ell}\Gamma_{q}\left(\ell \alpha+\beta\right)}\frac{\prod_{i=1}^{r}\left(a_{i};q\right)_{\ell}}{\prod_{j=1}^{s}\left(b_{j};q\right)_{\ell}}\frac{\Gamma_{q}\left(\nu+\ell\right)}{\Gamma_{q}\left(\nu+\mu+\ell\right)}z^{-\ell}q^{\mu\ell}\\
&=\frac{\Gamma_{q}(\nu)}{\Gamma_{q}\left(\nu+\mu\right)}\sum_{\ell=0}^{\infty}\frac{\left(-1\right)^{(1+s-r)(\ell)}q^{(1+s-r)\binom{\ell}{2}}}{\left(q;q\right)_{\ell}\Gamma_{q}\left(\ell \alpha+\beta\right)}\frac{\prod_{i=1}^{r}\left(a_{i};q\right)_{\ell}}{\prod_{j=1}^{s}\left(b_{j};q\right)_{\ell}}\frac{\left(q^{\nu};q\right)_{\ell}}{\left(q^{\nu+\mu};q\right)_{\ell}}z^{-\ell}q^{\mu\ell}\\
&=\frac{\Gamma_{q}\left(\nu\right)}{\Gamma_{q}\left(\nu+\mu\right)}\,_{r+1}R_{s+1,q}\left(a_{1},a_{2},\ldots,a_{r},q^{\nu};b_{1},b_{2},\ldots,b_{s},q^{\nu+\mu};\alpha,\beta;\frac{1}{z}q^{\mu}\right).
\end{split}
\end{equation*}
Similarly, we obtain result \eqref{6.22}.
\end{proof}
\begin{theorem}
The $q$-Weyl fractional $q$-integral and $q$-derivative operators for the $\,_{r}R_{s,q}$-function hold, defined as
\begin{equation}  \label{6.23}
\begin{split}
&W_{q}^{\nu,\mu}\left[\,_{r}R_{s,q}\left(a_{1},a_{2},\ldots,a_{r};b_{1},b_{2},\ldots,b_{s};\alpha,\beta;\frac{1}{z}\right)\right]\\
&=\frac{\Gamma_{q}(\nu)}{\Gamma_{q}\left(\nu+\mu\right)}\,_{r+1}R_{s+1,q}\left(a_{1},a_{2},\ldots,a_{r},q^{\nu};b_{1},b_{2},\ldots,b_{s},q^{\nu+\mu};\alpha,\beta;\frac{1}{z}q^{\mu}\right),
\end{split}
\end{equation}
and
\begin{equation} \label{6.24}
\begin{split}
&D_{q}^{\nu,\mu}\left[\,_{r}R_{s,q}\left(a_{1},a_{2},\ldots,a_{r};b_{1},b_{2},\ldots,b_{s};\alpha,\beta;z\right)\right]\\
&=\frac{\Gamma_{q}\left(\nu+\mu\right)}{\Gamma_{q}(\nu)}\,_{r+1}R_{s+1,q}\left(a_{1},a_{2},\ldots,a_{r},q^{\nu+\mu};b_{1},b_{2},\ldots,b_{s},q^{\nu};\alpha,\beta;\frac{1}{z}q^{-\mu}\right).
\end{split}
\end{equation}
\end{theorem}

\begin{proof}
 By virtue of the formulas \eqref{1.25} and \eqref{5.6}, we have
\begin{equation*}
\begin{split}
&W_{q}^{\nu,\mu}\left[\,_{r}R_{s,q}\left(a_{1},a_{2},\ldots,a_{r};b_{1},b_{2},\ldots,b_{s};\alpha,\beta;\frac{1}{z}\right)\right]\\
&=\sum_{\ell=0}^{\infty}\frac{c^{\ell}(-1)^{(1+s-r)(\ell)}q^{(1+s-r)\binom{\ell}{2}}}{\left(q;q\right)_{\ell}\Gamma_{q}\left(\ell \alpha+\beta\right)}\frac{\prod_{i=1}^{r}\left(a_{i};q\right)_{\ell}}{\prod_{j=1}^{s}\left(b_{j};q\right)_{\ell}}\,W_{q}^{\nu,\mu}z^{-\ell}\\
&=\sum_{\ell=0}^{\infty}\frac{(-1)^{(1+s-r)(\ell)}q^{(1+s-r)\binom{\ell}{2}}}{\left(q;q\right)_{\ell}\Gamma_{q}\left(\ell \alpha+\beta\right)}\frac{\prod_{i=1}^{r}\left(a_{i};q\right)_{\ell}}{\prod_{j=1}^{s}\left(b_{j};q\right)_{\ell}}\frac{\Gamma_{q}\left(\nu+\ell\right)}{\Gamma_{q}\left(\nu+\mu+\ell\right)}z^{-\ell}q^{\mu\ell}\\
&=\frac{\Gamma_{q}(\nu)}{\Gamma_{q}\left(\nu+\mu\right)}\sum_{\ell=0}^{\infty}\frac{(-1)^{(1+s-r)(\ell)}q^{\left(1+s-r\right)\binom{\ell}{2}}}{\left(q;q\right)_{\ell}\Gamma_{q}\left(\ell \alpha+\beta\right)}\frac{\prod_{i=1}^{r}\left(a_{i};q\right)_{\ell}}{\prod_{j=1}^{s}\left(b_{j};q\right)_{\ell}}\frac{\left(q^{\nu};q\right)_{\ell}}{\left(q^{\nu+\mu};q\right)_{\ell}}z^{-\ell}q^{\mu\ell}\\
&=\frac{\Gamma_{q}\left(\nu\right)}{\Gamma_{q}\left(\nu+\mu\right)}\,_{r+1}R_{s+1,q}\left(a_{1},a_{2},\ldots,a_{r},q^{\nu};b_{1},b_{2},\ldots,b_{s},q^{\nu+\mu};\alpha,\beta;\frac{1}{z}q^{\mu}\right).
\end{split}
\end{equation*}
Similarly, we obtain result \eqref{6.24}.
\end{proof}

\section{Special Cases}   \label{sec 7}
Here, we discuss the important special cases of the main results proved in the previous section, and we also explain the relevant aspects of the $\,_{r}R_{s,q}$ functions below.
\begin{theorem}
For $\Re(c)>\Re(b)>0$, the $\;_{2}R_{1,q}$ function satisfies the following $q$-integral representation:
\begin{equation}  \label{7.1}
	\begin{split}
	&\,_{2}R_{1,q}\left(a,b;c;\alpha,\beta;z\right)=\frac{\Gamma_{q}(c)}{\Gamma_{q}(b)\Gamma_{q}\left(c-b\right)}\int_{0}^{1}t^{b-1}\left(1-tq\right)^{c-b-1}E_{\alpha,\beta;q}^{a}\left(zt\right)\,\mathrm{d}_{q}t.
	\end{split}
	\end{equation}
	\end{theorem}
\begin{proof}
	For convenience, let the $\,_{2}R_{1,q}$ function denote the left-hand side of equation \eqref{7.1}. Then we have
\begin{equation*}
	\begin{split}
	&\,_{2}R_{1,q}\left(a,b;c;\alpha,\beta;z\right)=\sum_{\ell=0}^{\infty}\frac{z^{\ell}}{\left(q;q\right)_{\ell}\Gamma_{q}\left(\ell \alpha+\beta\right)}\frac{\left(a;q\right)_{\ell}\left(b;q\right)_{\ell}}{\left(c;q\right)_{\ell}}\\
	&=\frac{\Gamma_{q}(c)}{\Gamma_{q}(b)\Gamma_{q}\left(c-b\right)}\sum_{\ell=0}^{\infty}\frac{\left(a;q\right)_{\ell}z^{\ell}}{\left(q;q\right)_{\ell}\Gamma_{q}\left(\ell \alpha+\beta\right)}\int_{0}^{1}t^{b+\ell-1}\left(1-tq\right)^{c-b-1}\,\mathrm{d}_{q}t\\
&=\frac{\Gamma_{q}(c)}{\Gamma_{q}(b)\Gamma_{q}\left(c-b\right)}\int_{0}^{1}t^{b-1}\left(1-tq\right)^{c-b-1}\,E_{\alpha,\beta;q}^{a}\left(zt\right)\,\mathrm{d}_{q}t.
	\end{split}
	\end{equation*}
\end{proof}
\begin{remark}
Let $\alpha=1$, $\beta=1$ and $a_{r}=q$. Then \eqref{2.1} reduces to
\begin{equation*}
	\begin{split}
	\,_{r}R_{s,q}\left(a_{1},a_{2},\ldots,a_{r-1};q;b_{1},b_{2},\ldots,b_{s};1,1;z\right)=\,_{r-1}\Phi_{s}\left(a_{1},a_{2},\ldots,a_{r};b_{1},b_{2},\ldots,b_{s};q;\left(1-q\right)z\right),
	\end{split}
	\end{equation*}
where $\,_{r-1}\Phi_{s}$ is a generalized basic hypergeometric function with $r-1$ numerator and $s$ denominator parameters.
\end{remark}

\begin{remark}
Likewise, with $a_{1}=a$, $a_{2}=q$ $b_{1}=b$ and $\alpha=\beta=1$, the $R$--function reduces to basic confluent hypergeometric function, given as
\begin{equation*}
\,_{2}R_{1,q}\left(a,q;b;1,1;z\right)=\,_{1}\Phi_{1}\left(a;b;q;\left(1-q\right)z\right).
\end{equation*}
\end{remark}
\begin{remark}
	The classical $q$-Mittag-Leffler functions \cite{ab, abb} can also be represented as special cases of the $_{r}R_{s,q}$ function
\begin{equation}
\begin{split}
	&\,_{1}R_{0,q}\left(1;-;\alpha,1;z\right)=\sum_{\ell=0}^{\infty}\frac{\left(1;q\right)_{\ell}z^{\ell}}{\left(q;q\right)_{\ell}\Gamma_{q}\left(\ell \alpha+1\right)} = E_{\alpha,q}(z),
	\end{split}
	\end{equation}
\begin{equation}
	\begin{split}
	&\,_{1}R_{0,q}\left(1;-;\alpha,\beta;z\right)=\sum_{\ell=0}^{\infty}\frac{\left(1;q\right)_{\ell}z^{\ell}}{\left(q;q\right)_{\ell}\Gamma_{q}\left(\ell \alpha+\beta\right)} = E_{\alpha,\beta,q}(z),
	\end{split}
	\end{equation}
\begin{equation}
	\begin{split}
	&\,_{1}R_{0,q}\left(a;-;\alpha,\beta;z\right)=\sum_{\ell=0}^{\infty}\frac{\left(a;q\right)_{\ell}z^{\ell}}{\left(q;q\right)_{\ell}\Gamma_{q}\left(\ell \alpha+\beta\right)} = E_{\alpha,\beta,q}^{a}(z),
	\end{split}
	\end{equation}
and
\begin{equation}
	\begin{split}
	&\,_{2}R_{1,q}\left(a,1;c;\alpha,\beta;z\right) = \sum_{\ell=0}^{\infty}\frac{\left(a;q\right)_{\ell}\,z^{\ell}}{\left(q;q\right)_{\ell}\left(c;q\right)_{\ell}\Gamma_{q}\left(\ell \alpha+\beta\right)} = E_{\alpha,\beta,q}^{a,c}(z).
	\end{split}
	\end{equation}
\end{remark}
\begin{remark}
Further, the classical $q$-Wright function can be represented as a special case of the $\,_{0}R_{0,q}$ function with no numerator or denominator parameters, given by
\begin{equation}
\,_{0}R_{0,q}\left(-;-;\alpha,\beta;z\right) = \sum_{\ell=0}^{\infty}\frac{(-1)^{\ell}q^{\binom{\ell}{2}}z^{\ell}}{\left(q;q\right)_{\ell}\Gamma_{q}\left(\ell \alpha+\beta\right)}=\phi_{q}\left(\alpha;\beta;z\right).
	\end{equation}
\end{remark}
\begin{remark}
Let $a,b,c \in \mathbb{C}$, $\Re(a), \Re(b), \Re(c) >0$ and $\left|z\right|<1$, then we have
\begin{equation}
\,_{2}R_{1,q}\left(1,b;c;\alpha,1;z\right) = \frac{\Gamma_{q}(c)}{\Gamma_{q}(b)}\,_{2}\Psi_{2,q}\left(\left(1,1\right),\left(b,1\right);\left(c,1\right),\left(\alpha,1\right);z\right),
	\end{equation}
where$_{2}\Psi_{2,q}$ denotes the $q$-Wright function.
\end{remark}

\section{Conclusions and Further Directions}  \label{sec 8}
In the present paper, we investigate the properties, recurrence relations, and integral representations of the $_{r}R_{s,q}$ function. New results for the $q$-analogue of the $_{r}R_{s,q}$ function are derived by applying the $q$-Sumudu and $q$-Laplace transforms. The results presented here contribute new insights to the theory of $q$-fractional calculus and are likely to have useful applications across a wide range of problems in mathematics, statistics, and physical sciences. These results are an important addition to the theory of classical analysis, integral transforms, operational techniques, mathematical physics, fractional calculus, applied mathematics, and statistics. The findings in this paper are novel in the literature. We also provide several special cases, such as the $q$-hypergeometric, $q$-Mittag-Leffler, and $q$-Wright hypergeometric functions for the $_{r}R_{s,q}$ function, using the transform method. Additionally, we discuss their application to the $q$-Mittag-Leffler function and Euler-type integral representations, including various special cases, which will form the basis for future work.

\bigskip
\noindent \textbf{Declarations}\\
The authors have no conflicts of interest to declare that are relevant to the content of this article.

\noindent \textbf{Availability of data and materials}\\
No data were used for this study.

\noindent \textbf{Competing interests}\\
The authors declare that they have no competing interests.

\noindent \textbf{Funding}\\
No funding was received for conducting this study.

\end{document}